\documentclass[10pt,a4paper]{article}
\usepackage{amssymb, amsmath,amsfonts,amsthm}
\usepackage{graphicx}
\usepackage{color}
\usepackage[all]{xy}
\usepackage{psfrag}

\usepackage{combelow}

\theoremstyle{plain}
\newtheorem{theorem}{Theorem}[section]
\newtheorem{proposition}[theorem]{Proposition}
\newtheorem{corollary}[theorem]{Corollary}
\newtheorem*{claim*}{Claim}
\newtheorem{lemma}[theorem]{Lemma}

\theoremstyle{definition}

\newtheorem{notation}[theorem]{Notation}

\newcommand{\R}{\mathbb{R}}
\newcommand{\Q}{\mathbb{Q}}
\newcommand{\Z}{\mathbb{Z}}

\newcommand{\UT}{\mathrm{UT}}

\newcommand{\Cone}{\mathrm{C}1}
\newcommand{\Ctwo}{{\mathrm{C}2}}
\newcommand{\Cthree}{\mathrm{C}3}
\newcommand{\Cfour}{\mathrm{C}4}
\newcommand{\Cfive}{\mathrm{C}5}
\newcommand{\Rho}{P}

\DeclareMathOperator{\Aut}{Aut} 
\DeclareMathOperator{\ATF}{ATF} \DeclareMathOperator{\ITF}{ITF}

 \DeclareMathOperator{\im}{im}
 \DeclareMathOperator{\BS}{BS}

\newcommand{\pf}{\textbf{Proof: }}
\renewcommand{\qed}{\hfill{$\square$}}
\newcommand{\oag}{ordered abelian group}
\newcommand{\be}{\begin{enumerate}}
\newcommand{\ee}{\end{enumerate}}
\newcommand{\beq}{\begin{equation}}
\newcommand{\eeq}{\end{equation}}

\begingroup
    \catcode `\@ = 11
    \catcode `\~ = 13
    \catcode `\% = 12
    \protected\long\gdef\cmt@remove#1
    \ifdefined~
        \global\let\cmt@old~
    \else
        \global\let\cmt@old\relax
    \fi
    \protected\gdef~{\begingroup\catcode`
        \futurelet\next\cmt@}
    \protected\gdef\cmt@
      {\ifx
           \expandafter\cmt@remove
       \else
           \endgroup\expandafter\cmt@old
       \fi}
\endgroup

\renewcommand{\b}{\mathbf{b}}
\newcommand{\h}{\hspace*{0.4cm}}
\newcommand{\Y}{\mathrm{Y}}

\begin{document}
\bibliographystyle{plain}

\author{Shane O Rourke}
\title{A combination theorem for affine tree-free groups}\maketitle

\begin{abstract}
Let $\Lambda_0$ be an \oag. We show how an $\ATF(\Z\times\Lambda_0)$ group -- that is, a group admitting a free affine action without inversions on a $\Z\times\Lambda_0$-tree -- admits a natural graph of groups decomposition, where vertex groups inherit actions on $\Lambda_0$-trees. Using recent work of various authors,
it follows that a finitely generated group admitting a free affine action on a $\Z^n$-tree where no line has its orientation reversed is relatively hyperbolic with nilpotent parabolics, is locally quasiconvex, and has solvable word, conjugacy and isomorphism problems.

Conversely, given a graph of groups satisfying certain conditions, we show how an affine action of its fundamental group can be constructed. Specialising to the case of free affine actions, we obtain a large class of $\ATF(\Z\times\Lambda_0)$ groups that do not act freely by isometries on any $\Lambda_1$-tree.
We also give an example of a group that admits a free isometric action on a $\Z\times\Z$-tree but which is not residually nilpotent.\\

\end{abstract}
\section{Introduction}
A common theme in geometric theory is showing that if a
graph of groups is given whose vertex groups
belong to a particular class and the edge groups
are sufficiently well-behaved then the fundamental group of the graph
of groups also belongs to this class. Results of this sort have been established by Bestvina and Feighn \cite{comb-thm} in the case of hyperbolic groups, and by Dahmani \cite{Dahmani-combination}, Alibegovi{\'c} \cite{Alibegovic} and Bigdely and Wise \cite{Bigdely-Wise} in the case of relatively hyperbolic groups.

Our focus in this paper is on groups that admit actions of various sorts on $\Lambda_0$-trees. A group is $\ATF(\Lambda_0)$, or simply ATF, if it admits a free affine action (without inversions) on a $\Lambda_0$-tree. Similarly $\ITF(\Lambda_0)$ and ITF will be used to refer to groups that admit a free isometric action on a $\Lambda_0$-tree. We refer to the book \cite{Chiswell-book} for a detailed account of the fundamentals of $\Lambda_0$-trees and isometric actions thereon, and the survey paper \cite{KMS-long-survey} for an account of some of the fascinating recent developments in the theory of isometric actions on $\Lambda_0$-trees.

In \cite{Bass} Bass establishes fundamental results which show how
%
the isometric action of a group $\Gamma$ on a $\Z\times\Lambda_0$-tree gives rise to a graph of groups decomposition of $\Gamma$ together with isometric actions of the vertex groups on $\Lambda_0$-trees; conversely, there is a combination theorem to the effect that under certain natural conditions a collection of actions of vertex groups on $\Lambda_0$-trees can be extended to an action of the fundamental group (of a given graph of groups) on a $\Z\times\Lambda_0$-tree. (Here and throughout this paper we assume that a direct product or a direct power is endowed with the lexicographic order.) One can thus show that suitable combinations of $\ITF(\Lambda_0)$ groups are $\ITF(\Z\times\Lambda_0)$. In particular Martino and the author \cite{MOR} used Bass's results to show that certain groups admit free isometric actions on $\Z^n$-trees --- these include Liousse's examples of $\ATF(\R)$ groups described in \cite{Liousse}, as well as residually free surface groups. Bass's results have also been used by Kharlampovich, Miasnikov, Remeslennikov and Serbin \cite{KM-Zn} to give a detailed description of $\ITF(\Z^n)$ groups.

Affine actions were introduced in the case $\Lambda_0=\R$ by Liousse in \cite{Liousse}, and for general $\Lambda_0$ by the author in \cite{affine-paper}. We refer to the latter paper for the basic theory of affine actions on $\Lambda_0$-trees and properties of $\ATF(\Lambda_0)$ groups.

Briefly, if $X$ is a $\Lambda_0$-tree and $\beta_\sigma$ is an $o$-automorphism (order-preserving group automorphism) of $\Lambda_0$ and $\sigma$ is a permutation of $X$ with $d(\sigma x, \sigma y)=\beta_\sigma d(x,y)$ for all $x,y\in X$ then $\sigma$ is an \emph{affine automorphism} of $X$ (with dilation factor $\beta_\sigma$). One thus has the notion of a $\beta$-affine action of $\Gamma$ on a $\Lambda_0$-tree where $\beta:\Gamma\to\Aut^+(\Lambda_0)$ is a homomorphism and $\Aut^+(\Lambda_0)$ denotes the group of $o$-
automorphisms of $\Lambda_0$. While examples of groups that admit free isometric actions on $\Lambda_0$-trees for some $\Lambda_0$ include free groups, torsion-free abelian groups and fully residually free groups, examples of $\ATF$ groups include soluble Baumslag-Solitar groups $\BS(1,n)$ for $n\geq 1$, the wreath product $C_\infty\wr C_\infty$ and the Heisenberg group $\UT(3,\Z)$. Like their isometric counterparts, $\ATF$ groups are closed under free products and ultraproducts.

One key feature of Bass's combination theorem is that the translation lengths of (the embedded images of) an element of an edge group must match up. (Note that this requirement is also incorporated in the hypotheses of the main results of \cite{KM-Zn}.) In \cite{MOR} this was arranged in the cases considered by suitable adjustments to the metric in each of the $\Lambda_0$-trees on which the given end vertex groups acted. However such adjustments are not possible in general in the case of a circuit in the graph of groups: there are HNN extensions of free groups that cannot act freely by isometries on any $\Lambda$-tree for this reason. Examples of this phenomenon are furnished by
$$\Gamma=\langle F,t\ |\ tut^{-1}=v\rangle$$
where $F$ is a free group, and $u,v\in F$
generate distinct maximal cyclic subgroups of $F$ such that $u$ and $v$ cannot have the same translation length in any free action. An example of such $u$ and $v$ (attributed to Walter Parry) is described by Bass; a further example is $u=[x,y]$ and $v=[x^2,y^2]$, where $x$ and $y$ are non-commuting elements of $F$, as shown in \S\ref{one-rel-section} below.

Unfortunately the requirement that the translation lengths match up as just described cannot be neatly formulated in purely group-theoretic terms:
for the hyperbolic lengths $\ell(u)$ and $\ell(v)$ to be equal with respect to some free isometric action of a free group $F$ containing the elements $u$ and $v$ amounts to the existence of a free group $\bar{F}$ (on a basis $\bar{X}$) containing $F$ such that if $u_0$ and $v_0$ are cyclically reduced conjugates of $u$ and $v$, the word lengths (with respect to $\bar{X}$) of $u_0$ and $v_0$ are equal.

We establish affine analogues (Theorems~\ref{Bass3.5'} and~\ref{Bass3.8'} below) of Bass's theorems discussed above before specialising in \S3 to the case of free affine actions.
On the one hand, as in the isometric case, there is a natural hierarchical description of $\ATF(\Z^n)$ groups in terms of graphs of groups where the lowest ranked groups are free groups.
On the other hand, a notable difference between the affine case and the isometric case is that the homomorphism $\beta:\Gamma\to\Aut^+(\Lambda)$ affords an extra degree of freedom so that many groups that cannot admit free isometric actions do admit free affine actions on $\Z\times\Lambda_0$-trees.
This means that in many cases the awkward requirement concerning translation lengths of elements of the edge groups can be largely avoided giving somewhat cleaner corollaries than are possible in the isometric case. For example:

\begin{theorem}\label{intro-thm} Let $F$ be a free group and $u,v\in F$. The group $\Gamma=\langle F,t\ |\ tut^{-1}=v\rangle$ admits a free affine action on a $\Z\times\Q$-tree provided $u$ and $v$ are not proper powers and $v$ is not conjugate in $F$ to the inverse of $u$.
\end{theorem}

In fact, we will mainly consider \emph{essentially free} actions in \S3. This is a stronger condition than freeness, which is more robust in that it is preserved by two key constructions, namely the fulfilment of a $\Lambda_0$-tree, and the base change functor as applied to \emph{ample} embeddings of \oag s. In the isometric case a free action is automatically essentially free so these issues do not arise, but in our situation a discussion of ample embeddings and essentially free actions is necessary.

In \S\ref{rel-hyper-section} we use the graph of groups decomposition arising from a free affine action of a group $\Gamma$ on a $\Z^n$-tree as in Theorem~\ref{Bass3.5'} and the combination theorem of Bigdely and Wise \cite{Bigdely-Wise} to show that finitely generated $\ATF(\Z^n)$ groups are relatively hyperbolic where the parabolic groups are maximal nilpotent.
Now applying results of Farb \cite{Farb-rel-hyper}, Bumagin \cite{Bumagin}, and Dahmani and Touikan \cite{Dahmani-Touikan} we can deduce that finitely generated groups $\ATF(\Z^n)$ groups have solvable word, conjugacy and isomorphism problems respectively. Moreover, since finitely generated $\ATF(\Z^n)$ groups admit a so-called small hierarchy and nilpotent subgroups are Noetherian, by \cite[Theorem~D]{Bigdely-Wise} these groups are locally quasiconvex.

In \S\ref{one-rel-section} we consider the one-relator groups $\Gamma(m,n;r,s)=\langle x,y,t\ |\ t[x^m,y^n]t^{-1}=[x^r,y^s]\rangle$ and give necessary and sufficient conditions for these groups to admit free isometric or essentially free affine actions. As a by-product, we give an example of a finitely presented $\ITF(\Z^2)$ group that is not residually nilpotent.

Let us finally note that finitely presented ITF groups are automatic, and even bi-automatic (see \cite{KMS-long-survey}). By contrast, $\ATF(\Z^n)$ groups are not necessarily automatic: the Heisenberg group $\UT(3,\Z)$ is a well-known example of a non-automatic group, which is $\ATF(\Z^3)$.

On the other hand, the Heisenberg group is an example of a \emph{Cayley graph automatic} group, in the sense of Kharlampovich, Khoussainov and Myasnikov
(see \cite{KKM-CGA}). It would be interesting to know the relationship between $\ATF(\Z^n)$ groups and Cayley graph automatic groups, though we have not pursued this here.

I would like to thank Nicholas Touikan for kindly answering my questions about relatively hyperbolic groups.

\section{Affine actions on $\Lambda$-trees and graphs of groups}

A large part of the proofs of our results follow from arguments used in the isometric case in \cite{Bass}.
We will follow the notation of Bass \cite{Bass} as far as possible, and refer to this paper in case an argument used there can be routinely modified to apply in our situation.

We will use $\Lambda_0$ to refer to a fixed arbitrary \oag, and will generally use $\Lambda$ to refer to the \oag\ $\Z\times\Lambda_0$, equipped with the lexicographic order: thus $(m_1,\lambda_1)<(m_2,\lambda_2)$ if $m_1<m_2$, or $m_1=m_2$ and $\lambda_1<\lambda_2$. We will use $(X,d)$, or simply $X$, to refer to a $\Lambda$-tree. There is a natural projection $\Lambda\to\Z$, and a corresponding projection $p:X\to X^*$ ($x\mapsto x^*$) where $(X^*,d^*)$ is a $\Z$-tree: $X^*$ may be identified with the set of `balls of radius $0\times\Lambda_0 $'. Such a ball has the form $X(x^*)=p^{-1}(x^*)$. Of course $\Z$-trees may be viewed as trees in the usual sense: edges may be thought of as ordered pairs of vertices at distance 1 apart.
An action on $X$ induces an action on $X^*$ in a natural way. If the action on
$X^*$ is without inversions there is an associated quotient graph $Y^*$, a tree of representatives $T_0^*\subseteq X^*$ with $T_0^*$ isomorphic to a maximal subtree $T^*$ of $Y^*$, and a graph of groups $(\mathcal{G},Y^*,T^*)$ where vertex groups $\mathcal{G}(x^*)$ may be identified with certain vertex stabilisers $\Gamma_{x^*}$, and edge groups $\mathcal{G}(e)$ can be identified with certain edge stabilisers $\Gamma_e$. (The nature of these identifications depends on the choice of $T_0^*$ and a choice of group elements mapping one edge incident to $T_0^*$ to another such edge.)

We write $E(Y^*)$ for the set of edges of a graph $Y^*$, and $Y^*$ for the vertex set. Each edge $e$ is assumed to be oriented, with $\bar{e}$ denoting the oppositely oriented edge. Denote the origin of an edge $e$ by $\partial_0^{Y^*}e$, and the terminus of an edge by $\partial_1^{Y^*}e=\partial_0^{Y^*}\bar{e}$. We will occasionally drop the superscript $Y^*$ if it is clear from the context which graph we are considering. However, in the proofs of the results in this section, we will be considering both a tree $X^*$ and a quotient graph $Y^*$, and identifying vertices of $Y^*$ with those of a subtree of $X^*$, so in the interests of clarity we will usually write the superscript.

Two balls of radius $0\times\Lambda_0$ are adjacent if the corresponding vertices in $X^*$ are adjacent. If $x$ and $y$ belong to adjacent balls, then
$d(x,y)=(1,\lambda_0)$ (for some $\lambda_0$) and the segment $[x,y]$ may be expressed as a disjoint union $[x,\epsilon_{x})\cup [y,\epsilon_{y})$ where $[x,\epsilon_x)=\{z\in [x,y]: d(x,z)\in 0\times\Lambda_0\}$ and $[y,\epsilon_y)=\{z\in [x,y]: d(y,z)\in 0\times\Lambda_0\}$. The set $[x,\epsilon_x)$ is an $X(x^*)$-ray, and $\epsilon_x$ is an end of $X(x^*)$ of full $\Lambda_0$-type.

Let $X_0$ and $X_1$ be $\Lambda_0$-trees and suppose that $\epsilon_i$ is an end of $X_i$ of full $\Lambda_0$-type ($i=0,1$).
A \emph{bi-end map} (with associated ends $\epsilon_0$ and $\epsilon_1$) is a function $\Delta:X_0\times X_1\to\Lambda_0$ such that for $x_i\in X_i$ ($i=0,1$) the functions $x\mapsto\Delta(x,x_1)$ and $y\mapsto\Delta(x_0,y)$ are end maps (in the sense of \cite[1.4]{Bass}) towards $\epsilon_0$ and $\epsilon_1$ respectively. Note that conversely given two end maps $\delta_0$ and $\delta_1$ towards $\epsilon_0$ and $\epsilon_1$ respectively one can define a bi-end map by putting $\Delta(x,y)=\delta_0(x)+\delta_1(y)$.

Let $\delta_i$ and $\delta_i'$ be end maps towards $\epsilon_i$ ($i=0,1$).
By \cite[1.4]{Bass} there are constants $k_0$ and $k_1$ such that $\delta_i'=\delta_i+k_i$ ($i=0,1$).
Then, denoting the corresponding bi-end maps by $\Delta$ and $\Delta'$ respectively, we have $$\Delta'(x,y)=\delta_0'(x)+\delta_1'(y)=(\delta_0(x)+k_0)+(\delta_1(y)+k_1)=\Delta(x,y)+k_0+k_1.$$ Thus the bi-end maps $\Delta$ and $\Delta'$ are equal if and only if $k_1=-k_0$.

The reason for our interest in bi-end maps is that if $(X^*,d^*)$ is a $\Z$-tree, a $\Z\times\Lambda_0$-metric on $\coprod_{x^*\in X^*}X(x^*)$ which extends given $0\times\Lambda_0$-metrics on the $X(x^*)$, and whose projection onto $\Z$ induces $d^*$, is uniquely determined by a collection of bi-end maps $\Delta_e$ with associated ends $\epsilon_e$ and $\epsilon_{\bar{e}}$ joining the adjacent balls $X(\partial_0^{Y^*}e)$ and $X(\partial_0^{Y^*}\bar{e})$ ($e\in E(X^*)$) provided $\Delta_{\bar{e}}(y,x)=\Delta_e(x,y)$ for all $e$ (see \cite[2.2(b)]{Bass}). Conversely, the bi-end maps $\Delta_e$ are uniquely determined by $d$: for $x\in X(\partial_0 e)$ and $y\in X(\partial_0\bar{e})$ we have $d(x,y)=(1,-\Delta_e(x,y))$. We will sometimes speak of a bi-end map \emph{along $e$} for a bi-end map of the form $\Delta_e$
and write $\Delta_e s(x,y)$ for $\Delta_e(sx,sy)$;
that is, $s$ is understood to act diagonally.

In \cite[3.5]{Bass} Bass describes how an isometric action on a $\Z\times\Lambda_0$-tree can be `decomposed' as a graph of groups (arising from the action on the $\Z$-tree obtained by the projection $\Z\times\Lambda_0\to\Z$) equipped with actions of the vertex groups on $\Lambda_0$-trees satisfying certain natural compatibility conditions. We now wish to give an affine analogue of this result.
(Note that we include more data in our decomposition than Bass does. In particular, we will regard a choice of maximal subtree $T^*$ of $Y^*$ as part of the specification of a graph of groups.)\\

\emph{The given affine action}\\

Let $\Lambda_0$ be an \oag, $\Lambda=\Z\times\Lambda_0$ (lexicographically ordered), $\Gamma$ a group, and $\beta:\Gamma\to\Aut^+(\Lambda)$ a homomorphism so that $\beta_g(m,\lambda_0)=(m,\theta_g(\lambda_0)+m\mu_g)$ for some $\theta_g\in\Aut^+(\Lambda_0)$ and $\mu_g\in\Lambda_0$. Suppose that $\Gamma$ has a $\beta$-affine action on a $\Lambda$-tree $(X,d)$ for which the induced action on $X^*$ is without inversions, where $X^*$ is the
$\Z$-tree obtained from $X$ by identifying points whose distance apart is an
element of $0\times\Lambda_0$.\\

\emph{The graph of groups}\\

Then $\Gamma$ has an action on the $\Z$-tree $X^*$. Note that $X^*$ may be described in terms of the base change functor as $\Z\otimes_{\Lambda}X$. (See \cite[Theorem 8(3)]{affine-paper} for a discussion of the base change functor in the context of affine actions.)

Let $Y^*$ be the quotient graph of $X^*$ under the action of $\Gamma$ and
take a tree $T_0^*$ of representatives of $Y^*$ mod $\Gamma$ (see \cite[\S3.1]{Serre}). Thus $T_0^*$ is isomorphic to a maximal subtree $T^*$ of $Y^*$. We will later adopt the convention of identifying points of $T_0^*$ with their images under the projection map, and viewing $T^*$ as a subtree of $X^*$ as well as a (maximal) subtree of $Y^*$.

Extend $T_0^*$ to a subtree $S^*$ of $X^*$ such that the edge set $E(S^*)$ is mapped bijectively onto $E(Y^*)$ by the projection.
We will call such a subtree $S^*$ a tree of edge representatives.
Again we will often find it convenient not to distinguish between $E(S^*)$ and $E(Y^*)$.

For $e\in E(Y^*)$ choose $g_e\in\Gamma$ such that $g_e^{-1}\partial_0^{X^*}e\in T_0^*$; if $\partial_0^{X^*}e\in T_0^*$ we put $g_e=1$. (Note that for $e\in E(Y^*)$ either $g_e=1$ or $g_{\bar{e}}=1$.) Put $\mathcal{G}(x^*)=\Gamma_{x^*}$ for $x^*\in Y^*$ and $\mathcal{G}(e)=\Gamma_e$ for $e\in E(Y^*)$.

Let $\alpha_e:\mathcal{G}(e)\to\mathcal{G}(\partial_0^{Y^*}e)$ be given by $s\mapsto g_e^{-1}sg_e$. Then $(\mathcal{G},Y^*,T^*)$ is a graph of groups with fundamental group $\Gamma$.

For $e\in E(S^*)$, either $\partial_0^{X^*}e\in T_0^*$ or $\partial_0^{X^*}\bar{e}\in T_0^*$. If $\partial_0^{X^*}e\in T_0^*$ and $\partial_0^{X^*}\bar{e}\notin T_0^*$ we put $|e|=|\bar{e}|=e$; given the isomorphism between $T_0^*$ and $T^*$, this amounts to defining an orientation of $E(Y^*)\backslash E(T^*)$. One can write $\alpha_{|e|}(s)=s$ ($s\in\Gamma_e$) since $\alpha_{|e|}$ is a natural inclusion of $\Gamma_e$ in $\Gamma$.
If $\partial_0^{X^*}e,\partial_0^{X^*}\bar{e}\in T_0^*$ then $e$ is an edge of $T_0^*$, and $\alpha_e(s)$ and $\alpha_{\bar{e}}(s)$ are naturally identified in $\Gamma$. In this case we may still write $\alpha_{|e|}(s)$ for the common embedded image of $s\in\Gamma_e$ in $\Gamma$ under $\alpha_e$ and $\alpha_{\bar{e}}$.
\\

\emph{The roles of $\theta$ and $\mu$}\\

Note that $\beta_g(m,\lambda_0)=(m,\theta_g\lambda_0+m\mu_g)$ can be represented as a matrix equation $\beta_g\left(\begin{array}{c}\lambda_0\\ m\end{array}\right)=\left(\begin{array}{cc}\theta_g & \mu_g\\ 0 & 1\end{array}\right)\left(\begin{array}{c}\lambda_0\\ m\end{array}\right)=\left(\begin{array}{c}\theta_g\lambda_0+m\mu_g\\ m\end{array}\right)$.

The Bass-Serre relations for $\pi_1(\mathcal{G})$ give $\beta_{g_e}\beta_{\alpha_e(s)}\beta_{g_e}^{-1}=\beta_{g_{\bar{e}}}\beta_{\alpha_{\bar{e}}(s)}\beta_{g_{\bar{e}}}^{-1}$ and $\beta_{g_e}=1$ for $g_e=1$, which translate into the following conditions for $\theta_{g_e}$, $\theta^{x^*}$, $\mu_{g_e}$ and $\mu^{x^*}$. (Here $\theta^{x^*}$ and $\mu^{x^*}$ denote the respective restrictions of $\theta$ and $\mu$ to $\mathcal{G}(x^*)$, $x^*\in Y^*$.)

$$\begin{array}{rll}\theta_{g_e}\theta_{\alpha_e(g)}^{\partial_0e}\theta_{g_e}^{-1}
&=\theta_{g_{\bar{e}}}\theta_{\alpha_{\bar{e}(g)}}^{\partial_0\bar{e}}\theta_{g_{\bar{e}}}^{-1}, & s\in\mathcal{G}(e), e\in E(Y^*)\\
\theta_{g_e}&=1, &g_e=1\\
\mu_{g_e}+\theta_{g_e}\mu_{\alpha_e(s)}^{\partial_0e}-\theta_{g_e}\theta_{\alpha_e(s)}^{\partial_0e}\theta_{g_e}^{-1}\mu_{g_e}&=\mu_{\alpha_{|e|}(s)}^{\partial_0|e|}, & s\in\mathcal{G}(e), e\in E(Y^*)\\ \mu_{g_e}&=0, & g_e=1\end{array}$$

Moreover, $\theta^{x^*}:\mathcal{G}(x^*)\to\Aut^+(\Lambda_0)$ is a homomorphism, and $\mu^{x^*}$, while not a homomorphism, satisfies $\mu_{st}^{x^*}=\mu_s^{x^*}+\theta_s\mu_t^{x^*}$ ($s,t\in\mathcal{G}(x^*), x^*\in Y^*$).

Observe that $\beta_g$ is determined by its effect on elements of $\Lambda$ of the form $(1,\lambda_0)$ ($\lambda_0\in\Lambda_0$). In future we will generally write expressions such as $\beta_g:(1,\lambda_0)\mapsto(1,\theta_g\lambda_0+\mu_g)$ understanding that this defines $\beta_g$ on all of $\Lambda$.
\\

\emph{The $\Lambda_0$-trees $X(x^*)$ and the group actions thereon}\\

For $x^*\in T_0^*\subseteq X^*$ there is a $\Lambda_0$-tree $(X(x^*),d_{x^*})$ on which $\Gamma_{x^*}$ has a $\theta^{x^*}$-affine action where $\theta^{x^*}$ is the restriction of $\theta$ to $\Gamma_{x^*}$. Here $X(x^*)=p^{-1}(x^*)$
where $p :X\to X^*$ is the natural projection. Thus the fibres $p^{-1}(x^*)$ are precisely the balls of radius $0\times\Lambda_0$ in $X$
but we endow it with the $\Lambda_0$-metric $d_{x^*}(x,y)=\lambda_0$ where $d(x,y)=(0,\lambda_0)$.

For $x^*\in Y^*$ we thus obtain an associated $\Lambda_0$-tree equipped with an action of $\Gamma_{x^*}$ via the identification of $T_0^*$ with $T^*=Y^*$ (as vertex sets).
\\

\emph{The end maps}\\

For $e\in E(X^*)$ there is a bi-end map $\Delta_e:X(\partial_0^{X^*}e)\times X(\partial_0^{X^*}\bar{e})\to\Lambda_0$ with associated ends $\epsilon_e$ and $\epsilon_{\bar{e}}$, such that for $x,y\in X$ with $x^*$ adjacent in $X^*$ to $y^*$, we have $d(x,y)=(1,-\Delta_e(x,y))$.
These ends are of full $\Lambda_0$-type such that
$(\Gamma_{x^*})_{\epsilon_e}=\alpha_e\Gamma_e$.

Moreover, if $e,f\in E(Y^*)$ with $x^*=\partial_0^{Y^*}e=\partial_0^{Y^*}f$ then $\epsilon_e$ and $\epsilon_f$ are in distinct $\Gamma_{x^*}$-orbits unless $e=f$.

Now if $s$ is $\beta_s$-affine, where $\beta_s(1,\lambda_0)=(1,\theta_s\lambda_0+\mu_s)$, then for $x\in X(\partial_0e)$ and $y\in X(\partial_0\bar{e})$ we have

\begin{eqnarray*}d(sx,sy)&=&\beta_sd(x,y)\\ &=& \beta_s(1,-\Delta_e(x,y))\\ &=& (1,-\theta_s\Delta_e(x,y)+\mu_s)\end{eqnarray*} on the one hand, and $d(sx,sy)=(1,-\Delta_{se}(sx,sy))$
on the other.

Thus
$$\theta_s\cdot\Delta_e-\Delta_{se}\cdot s=\mu_s$$

For each pair of balls of radius $0\times\Lambda_0$ such as $X(x^*)$ and $X(y^*)$ adjacent via an edge $e$ say, one can choose end maps $\delta_e^{X^*}$ towards $\epsilon_e$ and $\delta_{\bar{e}}^{X^*}$ towards $\epsilon_{\bar{e}}$ such that $\Delta_e^{X^*}(x,y)=\delta_e^{X^*}+\delta_{\bar{e}}^{X^*}$. For $e\in E(S^*)$
we put $\delta_e^{Y^*}=\theta_{g_e}^{-1}\cdot\delta_e^{X^*}\cdot g_e$; thus if $\partial_0^{X^*}e\in T_0^*$ we have $\delta_e^{X^*}=\delta_e^{Y^*}$. Conversely, given end maps $\delta_e^{Y^*}$ for each $e\in E(Y^*)$ we can reverse these steps to define a metric on $\coprod X(x^*)$.
Note that unlike the bi-end maps $\Delta_e$, the end maps $\delta_e^{Y^*}$ are not uniquely determined by the metric on $X$: simultaneously replacing $\delta_e^{Y^*}$ by $\delta_e^{Y^*}+\theta_{g_e}^{-1}k$ and $\delta_{\bar{e}}^{Y^*}$ by $\delta_{\bar{e}}^{Y^*}-\theta_{g_{\bar{e}}}^{-1}k$ results in the same metric.

We can now rewrite the equation $\theta_s\cdot\Delta_e-\Delta_{se}\cdot s=\mu_s$ in terms of the end maps $\delta_e^{Y^*}$ as just described. Specialising to $s\in\mathcal{G}(e)$ gives the following equation.

$$\theta_{g_e}\left[\theta_{\alpha_e(s)}\delta_e^{Y^*}(x)-\delta_e^{Y^*}\cdot\alpha_e(s)(x)\right]
+\theta_{g_{\bar{e}}}\left[\theta_{\alpha_{\bar{e}}(s)}\delta_{\bar{e}}^{Y^*}(y)-\delta_{\bar{e}}^{Y^*}\cdot\alpha_{\bar{e}}(s)(y)
\right]=\mu_{\alpha_{|e|}(s)}\h s\in\mathcal{G}(e)$$

\emph{In short}\\

\begin{theorem}\label{Bass3.5'}
Let $\Lambda_0$ be an \oag\ and put
$\Lambda=\Z\times\Lambda_0$.

Suppose that\be
\item[(A1)] $\Gamma$ is a group;
\item[(A2)]$\beta:\Gamma\to\Aut^+(\Lambda)$
is a homomorphism where
$\beta_g(1,\lambda_0)=(1,\theta_g(\lambda_0)+\mu_g)$;
\item[(A3)] $X$ is a $\Lambda$-tree;
\item[(A4)] a
$\beta$-affine action of $\Gamma$ on $X$ is given for which the induced action on $X^*=\Z\otimes_{\Lambda}X$ is without inversions.
\ee

One can then obtain \be
\item[(S1)] a graph of groups $(\mathcal{G},Y^*,T^*)$ with vertex groups $\mathcal{G}(x^*)$,
edge groups $\mathcal{G}(e)$,
edge group embedding maps $\alpha_e$ ($e\in E(Y^*)$), elements $g_e\in\Gamma$ ($e\in E(Y^*)$), Bass-Serre tree $X^*$,
    and $T_0^*\subseteq S^*\subseteq X^*$, where $T_0^*$ is a tree of vertex representatives isomorphic to $T^*$, and $S^*$ is a tree of edge representatives.
The fundamental group $\pi_1(\mathcal{G})\cong\Gamma$ then has presentation $$\hspace{-2cm}\left\langle\ \mathcal{G}(x^*) (x^*\in Y^*),\ \ g_e (e\in E(Y^*)) \ |\ g_e\alpha_e(s)g_e^{-1}=g_{\bar{e}}\alpha_{\bar{e}}(s)g_{\bar{e}}^{-1}\ \ (e\in E(Y^*), s\in\mathcal{G}(e)),\ \
     g_e=1\ \ (\partial_0^{X^*}e\in T_0^*)\  \right\rangle.$$

Moreover, $g_e^{-1}\partial_0^{X^*}e\in T_0^*$ for $e\in E(S^*)$;

\item[(S2)] elements $\theta_{g_e}$ of $\Aut^+(\Lambda_0)$ ($e\in E(Y^*)$) and homomorphisms $\theta^{x^*}:\mathcal{G}(x^*)\to\Aut^+(\Lambda_0)$ ($x^*\in Y^*$), and
elements $\mu_{g_e}$ of $\Lambda_0$ ($e\in E(Y^*)$) and functions $\mu^{x^*}:\mathcal{G}(x^*)\to\Lambda_0$;
\item[(S3)] $\Lambda_0$-trees $X(x^*)$ for each $x^*\in Y^*$;
\item[(S4)] $\theta^{x^*}$-affine actions of $\mathcal{G}(x^*)$ on
$X(x^*)$ ($x^*\in Y^*$);

\item[(S5)]
ends $\epsilon_e$ of $X(\partial_0^{Y^*}e)$  ($e\in E(Y^*)$) and
end maps $\delta_e:X(\partial_0^{Y^*}e)\to\Lambda_0$ towards $\epsilon_e$ ($e\in E(Y^*)$).
\ee

These data satisfy the compatibility conditions
\be
\item[(C1)]  $\epsilon_e$ is of full $\Lambda_0$-type and $(\mathcal{G}(x^*))_{\epsilon_e}=\alpha_e\mathcal{G}(e)$ for $e\in E(Y^*)$ with $\partial_0^{Y^*}e=x^*$;
\item[(C2)]  If $e,f\in E(Y^*)$ with $x^*=\partial_0^{Y^*}e=\partial_0^{Y^*}f$ then $\epsilon_e$ and $\epsilon_f$ are in distinct $\mathcal{G}(x^*)$-orbits unless $e=f$;
\item[(C3)] $\theta_{g_e}\theta_{\alpha_e(s)}\theta_{g_e}^{-1}=\theta_{g_{\bar{e}}}\theta_{\alpha_{\bar{e}(s)}}\theta_{g_{\bar{e}}}^{-1}$ where $\theta_{\alpha_e(s)}=\theta_{\alpha_e(s)}^{\partial_0^{Y^*}e}$ and $\theta_{\alpha_{\bar{e}}(s)}=\theta_{\alpha_{\bar{e}(s)}}^{\partial_0^{Y^*}\bar{e}}$ ($s\in\mathcal{G}(e)$).
Moreover $\theta_{g_e}=1$ if $g_e=1$;
\item[(C4)] $\mu_{g_e}+\theta_{g_e}\mu_{\alpha_e(s)}-\theta_{g_e}\theta_{\alpha_e(s)}\theta_{g_e}^{-1}\mu_{g_e}=\mu_{\alpha_{|e|}(s)}$ ($s\in\mathcal{G}(e)$, $e\in E(Y^*)$), where $\mu_{\alpha_e(s)}=\mu_{\alpha_e(s)}^{\partial_0^{Y^*}e}$ ($s\in\mathcal{G}(e)$, $e\in E(Y^*)$). Also $\mu_{gh}=\mu_g+\theta_g\mu_h$ ($g,h\in\mathcal{G}(x^*)$, $x^*\in Y^*$) and $\mu_{g_e}=0$ if $g_e=1$;
\item[(C5)]  $\theta_{g_e}\left[\theta_{\alpha_e(s)}\delta_e^{Y^*}(x)-\delta_e^{Y^*}\cdot\alpha_e(s)(x)\right]
+\theta_{g_{\bar{e}}}\left[\theta_{\alpha_{\bar{e}}(s)}\delta_{\bar{e}}^{Y^*}(y)-\delta_{\bar{e}}^{Y^*}\cdot\alpha_{\bar{e}}(s)(y)
\right]=
\mu_{\alpha_{|e|}(s)}$\\ ($x\in X(\partial_0^{Y^*}e)$, $y\in X(\partial_0^{Y^*}\bar{e})$, $s\in\mathcal{G}(e)$, $e\in E(Y^*)$). \ee
\end{theorem}

\pf The data described in (S1)-(S5) all arise as described in the discussion preceding the theorem. That (C3)-(C5) are satisfied also follows from this discussion, while (C1) and (C2) follow just as in the isometric case; see \cite[3.5]{Bass}.
\qed

Given an action of a group $\Gamma$ on a $\Lambda$-tree, we will
call the data consisting of (S1)-(S5) the \emph{signature} of the
action.
Theorem~\ref{Bass3.5'} thus asserts that each
affine action on a $\Z\times\Lambda_0$-tree has a signature
satisfying (C1)-(C5).

Before tackling the converse to Theorem~\ref{Bass3.5'} we note that
to show a permutation of a
$\Z\times\Lambda_0$-tree is an affine automorphism, it suffices to
check points in the same ball or in adjacent
balls. More precisely:

\begin{lemma}\label{check-adj-balls}
Let $\Lambda_0$ be an \oag, put $\Lambda=\Z\times\Lambda_0$ (with
the lexicographic order) and let $X$ be a $\Lambda$-tree.
Let $s$ be a permutation of $X$ and $\beta_s\in\Aut^+(\Lambda)$ such
that
\begin{equation}\tag{$\ast$}\label{affine-adj} d(sx,sy)=\beta_s
d(x,y)\end{equation} whenever $d(x,y)=(m,\lambda_0)$ where $m\leq
1$. Then $s$ is a $\beta_s$-affine automorphism of $X$.
\end{lemma}

\pf
We claim that equation~(\ref{affine-adj}) is satisfied for all $x,y\in X$ and proceed by induction on $m=d^*(x^*,y^*)$. The required
assertion is given in the cases $m=0,1$, so assume that equation~(\ref{affine-adj}) is satisfied for all $l<m$.
Let $x=x_0,x_1,\ldots,x_m=y$ be points of $X$ such that
$[x,y]=[x_0,x_1,\ldots,x_m]$ and $d^*(x_{i-1}^*,x_{i}^*)=1$ for $1\leq i\leq m$. Thus $d(x,y)=\sum_{i=1}^m d(x_{i-1},x_i)$. Our hypothesis gives $d^*(sx_{i-1}^*,sx_i^*)=1$ for all $i$, so that $d^*(sx^*,sy^*)\leq m$; to show that we have equality here, let $u\neq v$ be points with $u^*=x_i^*=v^*$ and $[x_{i-1},u,v,x_{i+1}]$. Then $u\in[x_{i-1},v]$ and $d^*(x_{i-1}^*,u^*)=1=d^*(x_{i-1}^*,v^*)$ and $d^*(u^*,v^*)=0$, so by our hypothesis, \begin{eqnarray*}d(sx_{i-1},sv)&=&\beta_s d(x_{i-1},v)\\ &=&\beta_s\left(d(x_{i-1},u)+d(u,v)\right)\\ &=&\beta_s d(x_{i-1},u)+\beta_s d(u,v)\\ &=&d(sx_{i-1},su)+d(su,sv).\end{eqnarray*} That is, $su\in[sx_{i-1},sv]$. Similarly $sv\in[su,sx_{i+1}]$. Since $s$ is a permutation we have $su\neq sv$, so that $[sx_{i-1},su,sv,sx_{i+1}]$. Taking $u=x_i$ or $v=x_i$ we obtain $d(sx_{i-1},sx_{i+1})=d(sx_{i-1},sx_i)+d(sx_i,sx_{i+1})$. Therefore $[sx,sy]=[sx_0,sx_1,\ldots,sx_m]$, and \\ \begin{eqnarray*}d(sx,sy)&=&\sum_{i=1}^m d(sx_{i-1},sx_i)\\ &=&\sum_{i=1}^m\beta_s d(x_{i-1},x_i)\\ &=&\beta_s \sum_{i=1}^m d(x_{i-1},x_i)\\ &=&\beta_s d(x,y).\end{eqnarray*}
\qed

\begin{theorem}\label{Bass3.8'}

Let $\Lambda_0$ be an \oag\ and put $\Lambda=\Z\times\Lambda_0$.

Suppose that the following are given.

\be\item[(S1)]
a graph of groups $(\mathcal{G},Y^*,T^*)$ with Bass-Serre tree $X^*$, vertex groups $\mathcal{G}(x^*)$,
edge groups $\mathcal{G}(e)$,
edge group embedding maps $\alpha_e$ ($e\in E(Y^*)$),
    and $T_0^*\subseteq S^*\subseteq X^*$, where $T_0^*$ is a tree of vertex representatives isomorphic to $T^*$, and $S^*$ is a tree of edge representatives.
Thus $\pi_1(\mathcal{G})$ has presentation $$\hspace{-2cm}\left\langle\ \mathcal{G}({x^*}) (x^*\in Y^*),\ \ g_e (e\in E(Y^*)) \ | \ g_e\alpha_e(s)g_e^{-1}=g_{\bar{e}}\alpha_{\bar{e}}(s)g_{\bar{e}}^{-1}\ \ (e\in E(Y^*), s\in\mathcal{G}(e)),\ \
     g_e=1\ \ (\partial_0^{X^*}e\in T_0^*)\  \right\rangle$$
where the $g_e$ satisfy $g_e^{-1}\partial_0^{X^*}e\in T_0^*$ for $e\in E(S^*)$;
\item[(S2)] elements $\theta_{g_e}$ of $\Aut^+(\Lambda_0)$ ($e\in E(Y^*)$) and homomorphisms $\theta^{x^*}:\mathcal{G}(x^*)\to\Aut^+(\Lambda_0)$ ($x^*\in Y^*$), and
elements $\mu_{g_e}$ of $\Lambda_0$ ($e\in E(Y^*)$) and functions $\mu^{x^*}:\mathcal{G}(x^*)\to\Lambda_0$;
\item[(S3)]  $\Lambda_0$-trees $X(x^*)$ for $x^*\in Y^*$;
\item[(S4)]  $\theta^{x^*}$-affine actions of $\mathcal{G}(x^*)$ on $X(x^*)$ ($x^*\in Y^*$);
\item[(S5)]  ends $\epsilon_e$ of $X(\partial_0^{Y^*}e)$ and
end maps $\delta_e^{Y^*}:X(\partial_0^{Y^*}e)\to\Lambda_0$ towards $\epsilon_e$ ($e\in E(Y^*)$).
\ee

Assume that
\be\item[(C1)] $\epsilon_e$ is of full $\Lambda_0$-type and $\alpha_e\mathcal{G}(e)=\left(\mathcal{G}(\partial_0e)\right)_{\epsilon_e}$ ($e\in E(Y^*)$);
\item[(C2)]  if $\partial_0^{Y^*}e=x^*=\partial_0^{Y^*}f$ with $e\neq f$ then $\epsilon_e$ and $\epsilon_f$ lie in distinct $\Gamma_{x^*}$-orbits;
\item[(C3)]  $\theta_{g_e}\theta_{\alpha_e(s)}\theta_{g_e}^{-1}=\theta_{g_{\bar{e}}}\theta_{\alpha_{\bar{e}(s)}}\theta_{g_{\bar{e}}}^{-1}$, where $\theta_{\alpha_e(s)}=\theta_{\alpha_e(s)}^{\partial_0^{Y^*}e}$ and $\theta_{\alpha_{\bar{e}}(s)}=\theta_{\alpha_{\bar{e}}(s)}^{\partial_0^{Y^*}\bar{e}}$ for $e\in E(Y^*)$, $s\in\mathcal{G}(e)$. Moreover $\theta_{g_e}=1$ if $g_e=1$;
\item[(C4)] $\mu_{g_e}+\theta_{g_e}\mu_{\alpha_e(s)}-\theta_{g_e}\theta_{\alpha_e(s)}\theta_{g_e}^{-1}\mu_{g_e}=\mu_{\alpha_{|e|}(s)}$ ($s\in\mathcal{G}(e)$, $e\in E(Y^*)$), where $\mu_{\alpha_e(s)}=\mu_{\alpha_e(s)}^{\partial_0^{Y^*}e}$ ($s\in\mathcal{G}(e)$, $e\in E(Y^*)$). Also $\mu_{gh}=\mu_g+\theta_g\mu_h$ ($g,h\in\mathcal{G}(x^*)$, $x^*\in Y^*$) and $\mu_{g_e}=0$ if $g_e=1$;
\item[(C5)]      $\theta_{g_e}\left[\theta_{\alpha_e(s)}\delta_e^{Y^*}(x)-\delta_e^{Y^*}\cdot\alpha_e(s)(x)\right]
+\theta_{g_{\bar{e}}}\left[\theta_{\alpha_{\bar{e}}(s)}\delta_{\bar{e}}^{Y^*}(y)-\delta_{\bar{e}}^{Y^*}\cdot\alpha_{\bar{e}}(s)(y)
\right]=
\mu_{\alpha_{|e|}(s)}$\\ ($x\in X(\partial_0^{Y^*}e)$, $y\in X(\partial_0^{Y^*}\bar{e})$, $s\in\mathcal{G}(e)$, $e\in E(Y^*)$).
\ee

There exist \be\item[(A1)] a group $\Gamma$;
\item[(A2)] a homomorphism $\beta:\Gamma\to\Aut^+(\Lambda)$ given by $g\mapsto\beta_g$ where $\beta_g(1,\lambda_0)=(1,\theta_g\lambda_0+\mu_g)$;
\item[(A3)] a $\Lambda$-tree $X$;
\item[(A4)] a $\beta$-affine action of $\Gamma$ on $X$ with (S1)-(S5) as the signature.
\ee

\end{theorem}

\pf
\setcounter{equation}{0}
Much of the proof is similar to the proof given by Bass in \cite[3.8]{Bass} in the isometric case, and we will follow the argument and notation used there as far as possible. In particular we have given our equations labels (such as (B3)) to match up with equations (such as (3)) in Bass's proof.

If $\pi:X^*\to Y^*$ is the quotient map we will identify the vertex sets of $T_0^*$ and $Y^*$, and the edge sets $E(S^*)$ and $E(Y^*)$ via $\pi$; thus we will identify $T_0^*$ and $T^*$ as in (S1).
As in the isometric case we have $\partial_0^{X^*}e=g_e\partial_0^{Y^*}e$, giving
\begin{equation}\tag{B1}\label{Bass1} \partial_0^{X^*}(se)=sg_e\partial_0^{Y^*}(e)\h \mbox{for }e\in E(Y^*),\ s\in\pi_1(\mathcal{G}).
\end{equation}

We identify $t$ with its image $\alpha_e(t)$ in $\mathcal{G}(\partial_0e)$ in case $e\in E(S^*)$ and $\partial_0^{X^*}e\in T_0^*$. This gives
$$\alpha_e(t)=g_e^{-1}tg_e\h e\in E(S^*),\ t\in\mathcal{G}(e).$$

Let us now define the homomorphism $\beta$. Set $\Gamma=\pi_1(\mathcal{G})$. Then with respect to the action of $\Gamma$ on the Bass-Serre tree $X^*$, we have $\Gamma_{x^*}=\mathcal{G}(x^*)$ for $x^*\in Y^*$ and $\Gamma_e=\mathcal{G}(e)$ for $e\in E(S^*)$.

Putting $\beta_g=\beta_g^{x^*}:(1,\lambda_0)\mapsto(1,\theta_g^{x^*}\lambda_0+\mu_g^{x^*})$ ($g\in\mathcal{G}(x^*)$) and $\beta_{g_e}:(1,\lambda_0)\mapsto(1,\theta_{g_e}\lambda_0+\mu_{g_e})$, the compatibility conditions (C3) and (C4) ensure that the relations $\beta_{g_e}\beta_{\alpha_e(s)}\beta_{g_e}^{-1}=\beta_{g_{\bar{e}}}\beta_{\alpha_{\bar{e}}(s)}\beta_{g_{\bar{e}}}^{-1}$ ($s\in\mathcal{G}(e)$, $e\in E(Y^*)$) and $\beta_{g_e}=1$ ($\partial_0^{X^*}e\in T^*$) are satisfied. (We will henceforth omit the superscripts from $\mu_g$ and $\theta_g$.) The assignment $g\mapsto\beta_g$ therefore extends to a homomorphism
$\Gamma\to\Aut^+(\Lambda)$. It is easy to check that $\mu_{gh}=\mu_g+\theta_g\mu_h$ for all $g,h\in\Gamma$.

Now suppose that $\Lambda_0$-trees $(X(x^*),d_{x^*})$, acted upon by the vertex groups $\mathcal{G}(x^*)$ ($x^*\in Y^*$), and ends $\epsilon_e$ ($\partial_0^{Y^*}e=x^*$) are given as in (S3)-(S5).

We then have \begin{equation}\tag{B2}\label{Bass2}\mbox{If }e\in E(Y^*) \mbox{ and }x^*=\partial_0^{Y^*}(e) \mbox{ then } g_e^{-1}\Gamma_e g_e=\left(\Gamma_{x^*}\right)_{\epsilon_e}.\end{equation}

(Note that the second part of equation (2) as in \cite[3.8]{Bass} is not guaranteed in our situation, as the corresponding property is not assumed.)

We now define the set $X$ on which we will define the required $\Lambda$-metric. For $x^*\in Y^*$, we have an action of $\Gamma_{x^*}$ on $\Gamma\times X(x^*)$ via $t\cdot(s,x)=(st^{-1},tx)$. Denote the quotient arising from this action by $\Gamma\times_{\Gamma_{x^*}}X(x^*)$, and write $s\cdot x$ for the image of $(s,x)$ under the quotient map. We have an action of $\Gamma$ on $\Gamma\times_{\Gamma_{x^*}}X(x^*)$ by left multiplication. Now put $X(sx^*)=s\cdot X(x^*)$, and define a metric on $X(sx^*)$ by putting $d_{sx^*}(x_1,x_2)=\theta_s d_{x^*}(s^{-1}x_1,s^{-1}x_2)$. Note that if $sx^*=s_1x^*$ then $s^{-1}s_1\in\Gamma_{x^*}$,
so that $d_{x^*}(s^{-1}x_1, s^{-1}x_2)=\theta_{s^{-1}s_1}d_{x^*}(s_1^{-1}x_1,s_1^{-1}x_2)$, giving $\theta_s d_{x^*}(s^{-1}x_1, s^{-1}x_2)=\theta_{s_1}d_{x^*}(s_1^{-1}x_1,s_1^{-1}x_2)$. Thus the metric on $X(sx^*)$ depends only on $sx^*$, and not on $s$.

Now $$\Gamma\times_{\Gamma_{x^*}}X(x^*)=\coprod_{s\in\Gamma/\Gamma_{x^*}} X(sx^*),$$ and we put
$$X=\coprod_{x^*\in Y^*} \Gamma\times_{\Gamma_{x^*}}X(x^*),$$ and note that $\Gamma$ has an obvious action on $X$. Moreover $X=\coprod_{x^*\in X^*} X(x^*)$.

We define ends $\epsilon_{se}$, just as in \cite[3.8]{Bass}: \begin{equation}\tag{B3}\label{Bass3}\epsilon_{se}=sg_e\left(\epsilon_e^{Y^*}\right).\end{equation}
Using (\ref{Bass2}) above it is straightforward to check that this definition of $\epsilon_{se}$ depends only on $se$ not on the choice of $s$.

Since $0\times\Lambda_0$ is the maximal proper convex subgroup of $\Lambda$ (if $\Lambda_0\neq 0$), it is stabilised by $\beta_s$ for all $s\in\Gamma$. Thus $\epsilon_{se}$ is an end of full $\Lambda_0$-type in $X(\partial_0^{X^*}(se))$.

As in the isometric case one can show that
\begin{equation}\tag{B5}\label{Bass5}\mbox{If }te\mbox{ and }uf\mbox{ are
edges of }X^*\ (t,u\in\Gamma,\ e,f\in E(Y^*))\mbox{ and
}\epsilon_{te}=\epsilon_{uf} \mbox{ then } te=uf.\end{equation} The details are identical to the isometric case and we will omit them.
(In order to remain consistent with the
numbering of equations in \cite[3.8]{Bass}, we have also omitted equation (B4), since Bass's equation (4) does not adapt to the affine case.)

Taking the given end maps $\delta_e^{Y^*}$ and group elements $g_e$
($e\in E(Y^*)$), we define bi-end maps $\Delta_e^{X^*}$ ($e\in
E(X^*)$) as follows.

For $e\in E(S^*)$, we put $\delta_e^{X^*}=\theta_{g_e}\cdot\delta_e^{Y^*}\cdot g_e^{-1}$ and
define the bi-end map $\Delta_e^{X^*}(x,y)=\delta_e^{X^*}(x)+\delta_{\bar{e}}^{X^*}(y)$; it follows that $\Delta_{\bar{e}}(y,x)=\Delta_e(x,y)$.
We now
put
$$\Delta_{se}^{X^*}(x,y)=\theta_s\cdot\Delta_e(s^{-1}x,s^{-1}y)-\mu_s\h s\in\Gamma$$
Since every edge of $X^*$ has the form $se$ for some $e\in E(S^*)$ and $s\in\Gamma$, this defines $\Delta_e$ for all $e\in E(X^*)$.

\begin{claim*}
If $e\in E(S^*)$ and $se=s_1e$ then $\theta_s\cdot\Delta_e^{X^*}\cdot s^{-1}(x,y)-\mu_{s}=\theta_{s_1}\cdot\Delta_e^{X^*}\cdot s_1^{-1}(x,y)-\mu_{s_1}$.
\end{claim*}

\pf
Using the fact that $\alpha_e(s)=g_e^{-1}sg_e$, and putting $\bar{x}=g_e^{-1}s^{-1}x$ and $\bar{y}=g_{\bar{e}}^{-1}s^{-1}y$, we get

\begin{eqnarray*} \theta_s\cdot\Delta_e^{X^*}\cdot s^{-1}(x,y)-\Delta_e^{X^*}(x,y)
&=& \theta_s\cdot\delta_e^{X^*}\cdot s^{-1}(x)+\theta_s\cdot\delta_{\bar{e}}^{X^*}\cdot s^{-1}(y)-\delta_e^{X^*}(x)-\delta_{\bar{e}}^{X^*}(y)\\
&=& \theta_s\cdot \theta_{g_e}\cdot\delta_e^{Y^*}\cdot g_e^{-1}s^{-1}(x)+\theta_s\cdot \theta_{g_{\bar{e}}}\cdot\delta_{\bar{e}}^{Y^*}\cdot g_{\bar{e}}^{-1}s^{-1}(y)\\ & &\h
-\theta_{g_e}\cdot\delta_e^{Y^*}\cdot g_e^{-1}(x)-\theta_{g_{\bar{e}}}\cdot\delta_{\bar{e}}^{Y^*}\cdot g_{\bar{e}}^{-1}(y)\\
&=& \theta_{g_e}\left[\theta_{\alpha_e(s)}\cdot\delta_e^{Y^*}(\bar{x})-\delta_e^{Y^*}\cdot\alpha_e(s)(\bar{x})
\right]\\ & &\h +\theta_{g_{\bar{e}}}\left[
\theta_{\alpha_{\bar{e}(s)}}\cdot\delta_{\bar{e}}^{Y^*}(\bar{y})-\delta_{\bar{e}}^{Y^*}\cdot\alpha_{\bar{e}(s)}(\bar{y})
\right]\\
&=&\mu_{\alpha_{|e|}(s)}.
\end{eqnarray*}

Thus $\theta_s\cdot\Delta_e^{X^*}\cdot s^{-1}-\mu_{\alpha_{|e|}(s)}=\Delta_e^{X^*}=\theta_{s_1}\cdot\Delta_e^{X^*}\cdot s_1^{-1}-\mu_{\alpha_{|e|}(s_1)}$ for $s,s_1\in\Gamma_e$. If $\partial_0^{X^*}e\in T^*$ then $\alpha_{|e|}(s')=s'$ and the claim follows; otherwise it follows on replacing $e$ by $\bar{e}$ and noting that $\partial_0^{X^*}\bar{e}\in T^*$ and $\Delta_{\bar{e}}(y,x)=\Delta_e(x,y)$.
\medskip\\

It follows that $\Delta_{se}$ is well-defined.

We can now define a metric $d$ on $X$: for $x$ and $y$ belonging to a
common $\Lambda_0$-tree $X(x^*)=X(y^*)$ we put
$d(x,y)=(0,d_{x^*}(x,y))$, and for $x$ and $y$ belonging to distinct
$X(x^*)$ and $X(y^*)$ with $x^*=\partial_0^{X^*}e$ and $y^*=\partial_0^{X^*}{\bar{e}}$ we put
$$d(x,y)=(1,-\Delta_e(x,y)).$$ By \cite[2.2(b)]{Bass},
this suffices to specify the distance function on a
$\Z\times\Lambda_0$-tree. It remains to show that the action of
$\Gamma$ on $X$ is $\beta$-affine. Let $x,y\in X$ and $u\in\Gamma$.
Suppose first that $X(x^*)=X(y^*)$, $x^*,y^*\in T_0^*$ and $u$ stabilises the ball $X(x^*)$. Then
$u\in\Gamma_{x^*}=\mathcal{G}(x^*)$ and
$d_{x^*}(ux,uy)=\theta_u^{x^*}d_{x^*}(x,y)$ so that
$d(ux,uy)=(0,\theta_u^{x^*}d_{x^*}(x,y))=\beta_u(0,d_{x^*}(x,y))=\beta_ud(x,y)$.
If $u$ does not stabilise $X(x^*)$ then the definition of the
$\Lambda_0$-metric on $X(ux^*)=u\cdot X(x^*)$ ensures that $d(ux,uy)=\beta_u d(x,y)$. If $x^*,y^*\notin T_0^*$, then the required claim follows by taking $s\in\Gamma$ with $s^{-1}x,s^{-1}y\in T_0^*$. Then $d(ux,uy)=d(us(s^{-1}x),us(s^{-1}y))=\beta_{us}d(s^{-1}x,s^{-1}y)=\beta_u d(x,y)$.

Now suppose that $x$ and $y$ belong to adjacent balls:
$\partial_0^{X^*}e=x^*\neq y^*=\partial_0^{X^*}\bar{e}$, so that
$ux^*$ and $uy^*$ are the (distinct) endpoints of the edge $ue$ where $e\in E(S^*)$, and
$ux\in X(ux^*)$ and $uy\in X(uy^*)$. Using the definitions of $d$,
$\Delta_{ue}$ and $\delta_e^{X^*}$ we have
$$\begin{array}{rl}d(ux,uy)&=(1,-\Delta_{ue}^{X^*}(ux,uy))\\
&=(1,-\theta_u\Delta_e^{X^*}\cdot u^{-1}(ux,uy)+\mu_u)\\
&=\beta_u(1,-\Delta_e^{X^*}(x,y))\\ &=\beta_u d(x,y).
\end{array}$$

For $e\notin E(S^*)$ one can easily establish the same conclusion by considering $s\in\Gamma$ with $s^{-1}e\in E(S^*)$.
It follows from
Lemma~\ref{check-adj-balls} that $u$ is $\beta_u$-affine.\qed

Let us now consider the special case of Theorem~\ref{Bass3.8'} where the graph $Y^*$ has a single edge $e$. If $e$ has distinct endpoints $\partial_0e=x_0^*$ and $\partial_0\bar{e}=x_1^*$, and we put $\Gamma=\pi_1(\mathcal{G})$, then $\Gamma$ is an amalgamated free product $\Gamma_0\ast_{A_0=A_1}\Gamma_1$, where $A=\mathcal{G}(e)$, $\Gamma_i=\mathcal{G}(x_i^*)$ ($i=0,1$),
and $\alpha_e(a)=a_0$ and $\alpha_{\bar{e}(a)}=a_1$ are the edge group embeddings. Suppose that a $\theta^{x_i^*}$-affine action of $\Gamma_i$ on a $\Lambda_0$-tree $X(x_i^*)$ is given (where $\theta^{x_i^*}$ is a homomorphism), that $A_i$ is
the stabiliser of an end $\epsilon_i$ of $X(x_i^*)$ of full $\Lambda_0$-type ($i=0,1$). Suppose that end maps $\delta_e$ and $\delta_{\bar{e}}$ are given towards $\epsilon_0$ and $\epsilon_1$ respectively.
Assume that
$\theta^{x_0^*}_{a_0}=\theta^{x_1^*}_{a_1}$. Suppose also that $\mu_g\in\Lambda_0$ are given for $g\in \Gamma_i$ such that $\mu_{a_0}=\mu_{a_1}$ for $a\in A$ and $\mu_{gh}=\mu_g+\theta_g\mu_h$ for $g,h\in \Gamma_i$ ($i=0,1$).
Then, identifying $A$ with $A_0$ and $A_1$, there is a common extension $\theta:\Gamma\to\Aut^+(\Lambda_0)$ of the $\theta^{x_i^*}$ ($i=0,1$).
Moreover letting $\beta:\Gamma\to\Aut^+(\Lambda)$ be the homomorphism $\beta_g(1,\lambda_0)=(1,\theta_g(\lambda_0)+\mu_g)$, Theorem~\ref{Bass3.8'} now guarantees a $\beta$-affine action of $\Gamma$ on a $\Z\times\Lambda_0$-tree provided $$[\theta_a\delta_e(x)-\delta_e(ax)]+[\theta_a\delta_{\bar{e}}(y)-\delta_{\bar{e}}(ay)]=\mu_a\mbox{ for all }a\in A.
$$

Now suppose that $Y^*$ has one edge $e$ and a single vertex $x^*$. Then the fundamental group $\Gamma=\pi_1(\mathcal{G})$ is an HNN extension $\Gamma\cong\langle \Gamma_0,g_e\ |\ g_e\alpha_e(a)g_e^{-1}=a\ (a\in A)\rangle$ where $A=\mathcal{G}(e)$, and $\Gamma_0=\mathcal{G}(x^*)$. Suppose that a $\theta^{x^*}$-affine action of $\Gamma_0$ on a $\Lambda_0$-tree $X(x^*)$ is given (where $\theta^{x^*}$ is a homomorphism), that $A$ and $\alpha_e(A)$ are
the respective stabilisers of ends $\epsilon_0$ and $\epsilon_1$ of $X(x^*)$ of full $\Lambda_0$-type,
and that $\epsilon_0$ is in a distinct $\Gamma_0$-orbit from $\epsilon_1$. Suppose that end maps $\delta_e$ and $\delta_{\bar{e}}$ are given towards $\epsilon_0$ and $\epsilon_1$ respectively, that $\theta_{g_e}$, $\mu_{g_e}$ and $\mu_g$ ($g\in\Gamma_0$) are given, and assume that $\theta_{g_e}\theta^{x^*}_{\alpha_e(a)}\theta_{g_e}^{-1}=\theta^{x^*}_a$, $\mu_{g_e}+\theta_{g_e}\mu_{\alpha_e(a)}-\theta_{g_e}\theta_{\alpha_e(a)}\theta_{g_e}^{-1}\mu_{g_e}=\mu_a$ ($a\in A$), and $\mu_{gh}=\mu_g+\theta_g\mu_h$ ($g,h\in\Gamma_0$).

Let $\beta:\Gamma\to\Aut^+(\Lambda)$ be the homomorphism $\beta_g(1,\lambda_0)=(1,\theta_g(\lambda_0)+\mu_g)$.
Theorem~\ref{Bass3.8'} now guarantees a $\beta$-affine action of $\Gamma$ on a $\Z\times\Lambda_0$-tree provided
$$\theta_{g_e}[\theta_{\alpha_e(a)}\delta_e(x)-\delta_e\alpha_e(a)(x)]
+[\theta_a\delta_{\bar{e}}(y)-\delta_{\bar{e}}(ay)]=\mu_a\mbox{  for all }a\in A.
$$

\section{Free affine actions on $\Z\times\Lambda_0$-trees}

For us, an action of a group is \emph{free} if the stabiliser of every point and of every closed segment is trivial. In other words, we assume that a free action is without inversions.
We are particularly interested in the case where
our action is \emph{rigid}, that is, no segment is mapped properly into itself by any group element: this is because rigid automorphisms behave much like metric automorphisms.
We observe that if in the situation of Theorem~\ref{Bass3.5'} the action of $\Gamma$ is free then the induced action of $\Gamma_{x^*}=\mathcal{G}(x^*)$ on $X(x^*)$ is free. Conversely, if in the situation of Theorem~\ref{Bass3.8'} all vertex group actions are free, then so is the action of $\pi_1(\mathcal{G})$.

We will write $1$ for the identity automorphism of an \oag, and by extension, expressions such as $1-\theta_g$ for the map $\lambda_0\mapsto\lambda_0-\theta_g\lambda_0$. Recall that for an affine automorphism $\sigma$ of a $\Lambda_0$-tree we have $x\in A_\sigma$ if and only if $x\in[\sigma^{-1}x,\sigma x]$, which is in turn equivalent to the equation $d(\sigma^{-1}x,x)+d(x,\sigma x)=d(\sigma^{-1}x,\sigma x)$

\begin{proposition}\label{exist-free-action}
Let $\Lambda_0 $ be an \oag, $F$ a free group and $\theta:F\to\Aut^+(\Lambda_0 )$ a homomorphism. There exists a free $\theta$-affine action of $F$ on a $\Lambda_0 $-tree.
\end{proposition}

\pf
Let $\theta_*\in\Aut^+(\Lambda_0 )$ and suppose initially that $\theta_*$ is not the identity automorphism. Choose a positive $\lambda\in\Lambda_0 $ for which $\theta_*\lambda\neq\lambda$; without loss of generality, we will suppose that $\theta_*\lambda>\lambda$. Let $T_*$ be the subtree of $\Lambda_0 $ spanned by $\{\theta_*^k\lambda: k\in\Z\}$. For $\lambda_1\in T_*$, we have $\theta_*^m\lambda<\lambda_1<\theta_*^k\lambda$ for some distinct integers $m$ and $k$, giving $\theta_*^{-k}\lambda_1<\lambda<\theta_*^{-m}\lambda_1$, whence $\theta_*^{m-k}\lambda_1\neq\lambda_1$. Thus $\theta_*\lambda_1\neq\lambda_1$. Clearly $\theta_*$ preserves the order on $T_*$ inherited from $\Lambda_0 $, so that $\theta_*$ is a hyperbolic $\theta_*$-affine automorphism of $T_*$.

If $\theta_*$ is the identity automorphism, then for $\lambda_0>0$, the translation $\lambda\mapsto\lambda+\lambda_0$ is a $\theta_*$-affine (i.e. isometric) hyperbolic automorphism of $T=\Lambda_0 $.

Thus, for $\theta_x\in\Aut^+(\Lambda_0 )$ there is a free $\theta_x$-affine action of the infinite cyclic group $\langle x\rangle$ on a $\Lambda_0 $-tree. Viewing $F$ as the free product of cyclic groups, and applying \cite[Theorem 34]{affine-paper} we obtain a free $\theta$-affine action of $F$ on a $\Lambda_0 $-tree.
\qed

\subsection{Essentially free actions}\label{ess-free-actions}

We propose to describe a range of situations where a combination of ATF groups is ATF. It turns out that by making suitable modifications to the actions of the vertex groups, we can often arrange for the compatibility conditions (C1)-(C5) to be satisfied, even if the original actions do not satisfy them. For example, to apply Theorem 2.5, one needs the ends $\epsilon_e$ joining adjacent balls of radius $0\times\Lambda_0$ to be of full $\Lambda_0$-type. Using Bass's construction of the $\Lambda_0$-fulfilment of a $\Lambda_0$-tree (see \cite[Appendix E1]{Bass}), one can modify
the given $\Lambda_0$-trees so that all such ends are of full $\Lambda_0$-type.
Moreover an isometric automorphism $\sigma$ of a $\Lambda_0$-tree $X$ has a natural extension $\bar{\sigma}$ to its $\Lambda_0$-fulfilment $\bar{X}$, and $\bar{\sigma}$ is hyperbolic if $\sigma$ is. Thus a free (isometric) action on a $\Lambda_0$-tree extends to a free action on its $\Lambda_0$-fulfilment.

We now show that an affine action on a $\Lambda_0$-tree can be naturally extended to an action on its $\Lambda_0$-fulfilment.

\begin{lemma}\label{fulfilment}
Let $(X,d)$ be a $\Lambda_0$-tree, $\Gamma$ a group, $\theta:\Gamma\to\Aut^+(\Lambda_0)$ a homomorphism, and suppose that a $\theta$-affine action on $X$ is given. Let $(\bar{X},\bar{d})$ denote the $\Lambda_0$-fulfilment of $X$ (see \cite[E1]{Bass}), with associated embedding $\phi:X\to\bar{X}$. There exists a $\theta$-affine action of $\Gamma$ on $\bar{X}$ such that $\phi(gx)=g\phi(x)$ for $x\in X$ and $g\in\Gamma$.

In addition, the ends $\epsilon$ of $X$ are in canonical bijective correspondence with the ends $\epsilon_{\wedge}$ of $\bar{X}$ and the stabilisers satisfy $\Gamma_{\epsilon_{\wedge}}=\Gamma_{\epsilon}$. If $\delta:\bar{X}\to\Lambda$ is an end map towards $\epsilon_{\wedge}$ then $\delta$ restricted to (the embedded image of) $X$ (in $\bar{X}$) is an end map towards $\epsilon$ and for $s\in\Gamma_{\epsilon}$ with $\theta_s=1$, one has $\tau_{\epsilon}=\tau_{\epsilon_\wedge}$.
\end{lemma}

\pf
Fix $g\in\Gamma$, and
take $X=(X,d)$, $T=(\bar{X},\bar{d})$, $Y=(\bar{X},\theta_g^{-1}\bar{d})$ and $\psi=\psi_g=\phi\cdot g$. Then $\psi$ is an isometric embedding; moreover $g$ and $\phi$ are uniquely extending (see \cite[E1.6]{Bass}), and therefore, so is $\psi$. By \cite[E1.6]{Bass} there exists a unique surjective isometry $\psi'=\psi'_g:(\bar{X},\bar{d})\to(\bar{X},\theta_g^{-1}\bar{d})$ such that $\psi'\cdot\phi=\psi$. Moreover, uniqueness of $\psi'$ guarantees that $\psi'_{gh}=\psi'_g\psi'_h$. We therefore have an action of $\Gamma$ on $\bar{X}$ via $g\bar{x}=\psi'_g(\bar{x})$. Since $\psi'_g$ is an isometry with respect to the given metrics, it follows that $\psi'_g$ is a $\theta_g$-affine automorphism of $(\bar{X},\bar{d})$. Moreover since $\phi\cdot g=\psi_g=\psi'_g\cdot\phi$ we have $\phi(gx)=g\phi(x)$ for all $x\in X$.

The final assertions can be proven as in the isometric case: see \cite[E1.8]{Bass}.
\qed

However, given a hyperbolic affine automorphism of $X$, the induced automorphism of the $\Lambda_0$-fulfilment is not necessarily hyperbolic.
For example, if $\sigma$ is the hyperbolic affine automorphism $x\mapsto 2x$ of $X=(0,\infty)_{\R}$ then $\bar{\sigma}$ fixes $0\in\R=\bar{X}$.

We record also the following observation.

\begin{proposition}\label{o-action}
\be\item Let $G$ be a right-orderable (or equivalently, left-orderable) group. Then $G$ has a free affine action on a linear $\Lambda_1$ for some \oag\ $\Lambda_1$. If $G$ is orderable, the action is rigid.
\item If $G$ has a free rigid orientation-preserving action on a linearly ordered set then $G$ is orderable.
\ee
\end{proposition}

\pf
(a) Fix a compatible left order on $G$, and let $\Lambda_1$ be the subgroup of $\Z^G$ consisting of those elements $(n_g)=(n_g)_{g\in G}$ with finite support. The order on $G$ endows $\Lambda_1$ with an order (namely the lexicographic order), making $\Lambda_1$ an \oag. The natural action of $G$ on $\Lambda_1$ given by $\gamma:(n_g)_{g\in G}\mapsto(n_{\gamma g})_{g\in G}$ is $\beta$-affine (where $\beta_\gamma=\gamma$) and the set $X$ of positive elements of $\Lambda_1$ is invariant under the action.

If the left order on $G$ is also a right order then for positive $\gamma\in G$ and $(n_g)_{g\in G}\in X$ with $g_0=\min\{g\in G:n_g\neq 0\}$, we have $n_{g_0}>0$ and $n_g=0$ for $g<g_0$, and the first non-zero entry of $\gamma(n_g)_{g\in G}$ is also $n_{g_0}$, which appears in the $\gamma^{-1}g_0$ position. The right-invariance of the order now gives $\gamma^{-1}g_0<g_0$. This shows that $\gamma>1$ implies $\gamma(n_g)_{g\in G}>(n_g)_{g\in G}$. The converse can be shown by reversing the inequalities in the argument just given. Since $X$ consists of positive elements, it follows that the action of $G$ on $X$ is rigid.

To see that the action is free, suppose that $\gamma(n_g)=(n_g)$. If $g_0=\min\{g\in G:n_g\neq 0\}$ then $\gamma^{-1} g_0=g_0$, since this gives the respective positions of the first non-zero entries. Thus $\gamma=1$.

(b) Let a linearly ordered set $L$ be given together with an action as described, and let $x_0\in L$. Declare $g<h$ if $gx_0<hx_0$. This defines a linear order on $G$ which is left-invariant. In fact this order is independent of the choice of $x_0$ since $gx_0>x_0$ and $gx_1\leq x_1$ would violate rigidity. Thus $gx_0<hx_0$ implies $g\gamma x_0<h\gamma x_0$; in other words the order is right-invariant.
\qed

We will consider a condition (essential freeness) on an affine action that ensures that the induced action on the $\Lambda_0$-fulfilment is free. (The action described in Proposition~\ref{o-action} is not essentially free since zero is fixed by elements of $G$.) In fact our condition is somewhat stronger. Consider the map $g:x\mapsto\frac{x+4}{4}$; this defines an order-preserving $\theta_*$-affine automorphism of the dyadic rationals which has no fixed point. Therefore it is hyperbolic, although its extension (under the base change functor --- see \cite[Theorem 8(3)]{affine-paper}) to $\R$ fails to be hyperbolic. However, $g$ is not \emph{essentially hyperbolic} in a sense that we will presently make precise.

Let $\Lambda_0$ be an \oag, $\theta_*\in\Aut^+(\Lambda_0)$, and
$\lambda_0\in\Lambda_0$. We will say that $\lambda_0$ is
\emph{$\theta_*$-tame} if
$[\lambda_1-\theta_*(\lambda_1)]\subset[\lambda_0]$ for all
$\lambda_1\in\Lambda_0$. Note that if this condition is satisfied for
$\lambda_1$ with $\lambda_0\in[\lambda_1]$, it follows that the
condition is satisfied for all $\lambda_1$. For if $[\lambda_1]\subset[\lambda_0]$ and $[\lambda_1-\theta_*(\lambda_1)]\supseteq [\lambda_0]$, then $[\theta_*(\lambda_1)]\supseteq[\lambda_0]$. Replacing $\lambda_1$ by $\theta_*(\lambda_1)$, the assumed condition gives $[\theta_*(\lambda_1)-\theta_*^2(\lambda_1)]\subset[\lambda_0]$. On the other hand we have $[\lambda_1]\subset[\lambda_0]\subseteq[\theta_*(\lambda_1)]$, giving $[\theta_*(\lambda_1)]\subset[\theta_*^2(\lambda_1)]$ and hence $[\theta_*(\lambda_1)-\theta_*^2(\lambda_1)]=[\theta_*^2(\lambda_1)]\supset[\lambda_0]$, a contradiction.

We will find it convenient to use the notation $a\ll b$, or equivalently $b\gg a$, if $ka<b$ for all $k\in\Z$. In particular, $b$ must be positive. We will write $[a]=[a]_{\Lambda_0}$ for the convex subgroup of $\Lambda_0$ generated by $a$. It is easy to check that \be\item if $a\ll b$ and $b\leq c$ then $a\ll c$;
\item if $a,b\ll c$ then $ka+lb\ll c$ for $k,l\in\Z$; \item if $a,b\ll c$ then $a\ll kb+lc$ for $k,l\in\Z$, $l>0$.\ee However $a\ll |b|,|c|$ does not imply $a\ll |b+c|$ as the case $c=-b$ shows. We will use properties (a)-(c) throughout this section without further comment.
Note also that $$[a]_{\Lambda}\subset[b]_{\Lambda}\Leftrightarrow ka<|b|\mbox{\ for all }k\in\Z\Leftrightarrow a\ll |b|\Leftrightarrow |a|\ll|b|.$$

\begin{lemma}\label{axis-embed}
Let $\sigma$ be an order-preserving $\theta_{\sigma}$-affine automorphism of a linear $\Lambda_0$-tree $L$, and let $\iota:L\to\Lambda_0$ be an isometric embedding. There exists $\nu_{\sigma}\in\Lambda_0$ such that $$\iota(\sigma x)=\theta_{\sigma}\iota(x)+\nu_{\sigma}\h \mbox{for all }x\in L.$$ Moreover,
$\nu_\sigma$ and $\nu_{\sigma^{-1}}$ have opposite signs.
\end{lemma}
\pf
Fix a compatible linear order $\leq$ on $L$ and let $x_0\in L$. For $x\in L$ such that $\iota(x)\geq \iota(x_0)$ we have $\iota(\sigma x)\geq \iota(\sigma x_0)$, whence
\begin{eqnarray*}\iota(\sigma x)-\iota(\sigma x_0)&=&d(\sigma x,\sigma x_0)\\ &=&\theta_{\sigma}d(x,x_0)\\ &=&\theta_{\sigma}(\iota(x)-\iota(x_0))\\ &=& \theta_{\sigma}\iota(x)-\theta_{\sigma}\iota(x_0);\end{eqnarray*} the equation $\iota(\sigma x)-\iota(\sigma x_0)=\theta_\sigma\iota(x)-\theta_\sigma\iota(x_0)$ is similarly established if $\iota(x)\leq\iota(x_0)$. Thus $\iota(\sigma x)-\theta_{\sigma}\iota(x)$ is independent of the choice of $x\in L$.
Thus there is a constant $\nu_{\sigma}$ such that $\iota(\sigma x)-\theta_{\sigma}\iota(x)=\nu_{\sigma}$.

To establish the final assertion, note that $\nu_{\sigma^{-1}}>0$ if and only if $\iota(\sigma^{-1}x')>\theta_{\sigma^{-1}}\iota(x')$ for all $x'\in A_{\sigma^{-1}}=A_{\sigma}$. If $\nu_{\sigma}$ and $\nu_{\sigma^{-1}}$ are both positive then
$\iota(\sigma x)>\theta_{\sigma}\iota(x)$ for all $x\in A_\sigma$. Putting $x'=\sigma x$ then yields a contradiction. The case where both $\nu_{\sigma^{-1}}$ and $\nu_{\sigma}$ are negative can be similarly dismissed.
\qed

Note that $\nu_{\sigma}$ depends on the choice of isometric embedding $\iota$. Replacing $\iota$ by $\pm\iota+k$ has the effect of replacing $\nu_{\sigma}$ by $\pm\nu_{\sigma}+(1-\theta_{\sigma})(k)$. This accounts for all possible values of $\nu_\sigma$ since all other isometric embeddings are of this form; this follows from \cite[Lemma 1.2.1]{Chiswell-book}.

A hyperbolic $\theta_\sigma$ automorphism $\sigma$ of a $\Lambda_0$-tree $X$ is \emph{essentially hyperbolic} if $\nu_\sigma$ is $\theta_\sigma$-tame; explicitly, if $[(1-\theta_\sigma)(\lambda_0)]_{\Lambda_0}\subset[\nu_\sigma]_{\Lambda_0}$ for all $\lambda_0\in\Lambda_0$. Thus $\sigma$ is essentially hyperbolic if $(1-\theta_\sigma)(\lambda_0)\ll|\nu_\sigma|$ for all $\lambda_0$, which, for $\lambda_1\in\Lambda_0$, is equivalent to the requirement that $(1-\theta_\sigma)(\lambda_0+\lambda_1)\ll|\nu_\sigma\pm(1-\theta_\sigma)(\lambda_1)|$ for all $\lambda_0$. It follows that the designation of $\sigma$ as essentially hyperbolic is independent of the choice of $\iota$ that determines the value of $\nu_\sigma$. Note that a hyperbolic isometry is automatically essentially hyperbolic.

We will write $\lambda'=\lambda+o(g)$ if $\lambda'-\lambda\in[\im(1-\theta_g)]$; thus, $g$ is essentially hyperbolic if and only if $\nu_g\neq 0+o(g)=o(g)$. (Since all groups under consideration are torsion-free, this should cause no confusion with the order of $g$.) An affine action is \emph{essentially free} if all non-trivial elements are essentially hyperbolic.

\begin{lemma}\label{ess-hyper-props}
\be\item Let $\sigma$ be a hyperbolic non-rigid affine automorphism of a $\Lambda_0 $-tree. There exist $x,y\in A_\sigma$ for which $\sigma[x,y]\subset[x,y]$ or $[x,y]\subset\sigma[x,y]$.
\item Essentially hyperbolic automorphisms are rigid.
\item A hyperbolic affine automorphism of an $\R$-tree is essentially hyperbolic if and only if it is an isometry.
\item If $\Gamma$ has an essentially free affine action on a $\Lambda_0$-tree $X_0$ then the induced affine action on its $\Lambda_0$-fulfilment is essentially free (and hence free).
\item Let $g,h\in G$ where $G$ has an essentially free affine action on a $\Lambda_0$-tree and suppose that $h^{-1}gh=g^{-1}$. Then $g=1$.
\ee
\end{lemma}

\pf
(a) Replacing $\sigma$ by its inverse if necessary, we may suppose that $x,y\in X$ with $\sigma[x,y]\subset[x,y]$. Replacing $\sigma$ by $\sigma^2$ we may assume that $[x,\sigma x,\sigma y,y]$; moreover $A_\sigma=A_{\sigma^2}$, by \cite[Theorem 14(10)]{affine-paper}. By \cite[Theorem 14(8)]{affine-paper} there exists a point $p\in[x,\sigma x]\cap A_\sigma$ and a point $q\in[\sigma y,y]$. Since $A_\sigma$ is a subtree, we therefore have $\sigma x,\sigma y\in A_\sigma$, and thus $x,y\in A_\sigma$.

(b) Suppose that $\sigma$ is a hyperbolic automorphism of $X$ but is not rigid. Then embedding $A_\sigma$ in $\Lambda$ as in Lemma~\ref{axis-embed}, there exist $x_0,x_1\in A_\sigma$ such that either $\theta_\sigma\iota(x_0)+\nu_\sigma=\iota\sigma(x_0)\leq\iota(x_0)<\iota(x_1)<\iota\sigma(x_1)=\theta_\sigma\iota(x_1)+\nu_\sigma$ or $\iota(x_0)\leq \theta_\sigma\iota(x_0)+\nu_\sigma<\theta_\sigma\iota(x_1)+\nu_\sigma<\iota(x_1)$. In the former case we then have $(1-\theta_\sigma)(\iota(x_1))<\nu_\sigma\leq(1-\theta_\sigma)(\iota(x_0))$, so that $[\im(1-\theta_\sigma)]$ cannot be properly contained in $[\nu_\sigma]$ for all $\lambda$; thus $\sigma$ cannot be essentially hyperbolic. The latter case can be handled similarly.

(c) If $(1-\theta_\sigma)(\lambda)\ll |\nu_\sigma|$ for all $\lambda\in\R$, then this forces $(1-\theta_\sigma)(\lambda)=0$; that is, $\theta_\sigma=1$. The converse is trivial.

(d)
Consider the induced $\theta$-affine action of $\Gamma$ on the $\Lambda_0$-fulfilment $\bar{X_0}$ of $X_0$. For non-trivial $g\in \Gamma$, consider $A_g\subseteq \bar{X_0}$  and let $\iota:A_g\to\Lambda_0$ be an injection; since $\bar{X_0}$ is full, $\iota$ is in fact surjective. If $gx=x$ for some $x\in A_g$, then $\iota(x)=\iota(gx)=\theta_g\iota(x)+\nu_g$, giving $(1-\theta_g)(\iota(x))=\nu_g$, which violates the essential hyperbolicity of $g$.

(e)
By part (b) an essentially free action is free and rigid. Now $A_g=A_{g^{-1}}=A_{h^{-1}gh}=h^{-1}A_g$, so that $A_g$ is $h$-invariant. If $g\neq 1$ then $g$ is hyperbolic and $A_g=A_h$ is a linear subtree stabilised by $\langle g,h\rangle$; moreover the orientation is preserved by $g$ and $h$. By Proposition~\ref{o-action}(b) $\langle g,h\rangle$ is orderable, which is impossible, since only the identity can be conjugate to its inverse in an orderable group.
\h \qed

Note that the converse of Lemma~\ref{ess-hyper-props}(b) is not valid: a simple example is afforded by $X=\{0\}\times\Z$, $\Lambda=\Z\times\Z=\bar{X}$, $\theta_\sigma(k,m)=(k,m+k)$, and $\sigma(k,m)=(k,m+k+1)$.

There are cases of $\theta_*\in\Aut^+(\Lambda_0)$ for which there is no essentially hyperbolic $\theta_*$-affine automorphism of any $\Lambda_0$-tree: consider $\Lambda_0=\omega(\Z^\Z)$, and $\theta_*:(n_k)_{k\in\Z}\mapsto(n_{k+1})_{k\in\Z}$. Then for any positive $(\nu_k)_{k\in\Z}\in \Lambda_0$ with least non-zero entry in the $k_0$th position we can take $n_{k_0}$ to be negative and all other $n_k$ to be zero. Then $(1-\theta_*)(n_k)>(\nu_k)$; thus no element of $\Lambda_0$ is $\theta_*$-tame.
It follows that Proposition~\ref{exist-free-action} cannot be strengthened to guarantee the existence of non-trivial essentially free actions on $\Lambda_0$-trees for all \oag s $\Lambda_0$.

In \cite{affine-paper} we showed that free affine actions of groups yield free affine actions of their free products and ultraproducts. We now show that essential freeness is similarly preserved by these constructions.

\begin{proposition}
\be\item For $x\in X$ and $u=\Y(g^{-1}x,x,gx)$, let $\b_x(g)=d(u,gu)$. We have \be\item $\nu_g=\epsilon\b_x(g)+o(g)$ where $\epsilon=\pm 1$;

\item $(1-\theta_{\gamma g\gamma^{-1}})=\theta_{\gamma}(1-\theta_g)\theta_{\gamma^{-1}}$;

\item $\b_{\gamma x}(\gamma g\gamma^{-1})=\theta_{\gamma}\b_x(g)$;

\item $\theta_\gamma o(g)=o(\gamma g\gamma^{-1})$;

\item $\nu_{\gamma g\gamma^{-1}}=\pm \theta_\gamma\nu_g+o(\gamma g\gamma^{-1})$;

\item If $g_1,\ldots,g_m\in G$ and $\overline{g_k}$ is the product $g_1\cdots g_k$, then $(1-\theta_{\overline{g_m}})=\sum_{k=1}^m\theta_{\overline{g_{k-1}}}(1-\theta_{g_k})$.
\ee

\item Suppose that $G_i$ has an essentially free $\theta^{(i)}$-affine action on a $\Lambda_0 $-tree for $i\in I$ (where $\theta^{(i)}:G_i\to\Aut^+(\Lambda_0 )$ is a homomorphism.) The induced action of the free product $G=\ast_{i\in I}G_i$ on the $\Lambda_0 $-tree $X$ (see \cite[Theorem 34]{affine-paper}) is essentially free.
\ee
\end{proposition}

\pf
(a)
(i) If $\iota:A_g\to\Lambda_0 $ is chosen so that $\iota(gu)\geq\iota(u)=0$, we have \begin{eqnarray*}\nu_g&=&\iota(gu)-\theta_g\iota(u)\\ &=&\iota(gu)-0\\ &=& \iota(gu)-\iota(u)\\ &=& \b_x(g),
\end{eqnarray*} and for any other isometric embedding $\iota'$, we have $\nu_g'=\epsilon\b_x(g)+(1-\theta_g)(k)= \epsilon\b_x(g)+o(g)$, where $\epsilon=\pm 1$.

Assertions (ii) and (iii) are routine; let us consider (iv). Note first that if $\lambda_2=\theta_\gamma\lambda_1$ then $(1-\theta_g)(\lambda_1)$ and $(1-\theta_{\gamma g\gamma^{-1}})(\lambda_2)=\theta_\gamma(1-\theta_g)(\lambda_1)$ have the same sign. Thus $\lambda_0<|(1-\theta_g)(\lambda_1)|$ precisely when $\theta_\gamma(\lambda_0)<|(1-\theta_{\gamma g\gamma^{-1}})(\lambda_2)|$. Therefore $\lambda_0=o(g)$ precisely when $\theta_\gamma(\lambda_0)=o(\gamma g\gamma^{-1})$.

(v)
From part (i) and (iii),
we have (for some $\epsilon,\epsilon'\in\{1,-1\}$) \begin{eqnarray*}\nu_{\gamma g\gamma^{-1}}&=&\epsilon\theta_\gamma\b_{\gamma^{-1}x}(g)+o(\gamma g\gamma^{-1})\\ &=&\epsilon\theta_\gamma(\epsilon'\nu_g+o(g))+o(\gamma g\gamma^{-1})\\ &=& \pm\theta_\gamma\nu_g+o(\gamma g\gamma^{-1})
\end{eqnarray*}

Part (vi) is a routine calculation.

(b) The cited theorem guarantees that the action of $G$ is free. We will use the notation $g\cdot h$ to denote the product $gh$ in the case where $g$ and $h$ belong to distinct free factors.
Let $g=g_1\cdot g_2\cdots g_k$. If $k=1$ then the required assertion $(1-\theta_g)(\lambda_0)\ll|\nu_g|$ is immediate from our hypotheses. So assume that $k\geq 2$.
We can write $g$ as a conjugate $\gamma (h_1\cdot h_2\cdots h_m)\gamma^{-1}$ where $h_m$ and $h_{1}$ belong to distinct free factors. Now $\nu_g=\pm\theta_\gamma\nu_{h_1\cdots h_m}+\theta_\gamma o(h_1\cdots h_m)$ by parts (v) and (iv), and $(1-\theta_g)(\lambda_0)=\theta_{\gamma}(1-\theta_{\overline{h_m}})(\theta_{\gamma^{-1}}\lambda_0)$, using part (ii). It therefore suffices to show that $(1-\theta_g)(\lambda_0)\ll|\nu_g|$ in the case where $g=g_1\cdot g_2\cdots g_k$, and $g_k$ and $g_1$ belong to distinct free factors.

As in \cite[Theorem 34]{affine-paper}, we choose a basepoint $x_i$ for each action of $G_i$ ($i\in I$), and let $x$ be the basepoint of the resulting action of $G$ on $X$ for which $L=L_x$. Then taking $g=g_1\cdots g_k$, we have $d(g^{-1}x,x)+d(x,gx)=\theta_{g^{-1}}(1+\theta_g)L(g)$, while $g^2=g_1\cdots g_k\cdot g_1\cdots g_k$, whence $L(g^2)=L(g)+\theta_g L(g)=(1+\theta_g)L(g)$ by \cite[Theorem 34]{affine-paper}. It follows that $d(g^{-1}x,x)+d(x,gx)=\theta_{g^{-1}}L(g^2)=d(g^{-1}x,gx)$; thus $x\in A_g$. Therefore $L(g)=\b_x(g)=\epsilon\nu_g+o(g)$. Moreover, our hypotheses ensure that $(1-\theta_{g_i})(\lambda_0)\ll|\nu_{g_i}|=\b_x(g_i)\leq L(g_i)$.
Since by definition $L$ satisfies $L(g_1\cdot g_2\cdots g_k)=\sum_{i=1}^k\theta_{\overline{g_{i-1}}}L(g_{i})$ for $k\geq 1$, we have \begin{eqnarray*}(1-\theta_g)(\lambda_0)&=&\sum_{i=1}^k \theta_{\overline{g_{i-1}}}(1-\theta_{g_i})(\lambda_0)\\  &\ll& \sum_{i=1}^k \theta_{\overline{g_{i-1}}}L(g_i)\\ &=& L(g)\\ &=&|\nu_g|+o(g)
\end{eqnarray*}
\qed

\begin{proposition}
If $G_i$ has an essentially free $\theta^{(i)}$-action on a $\Lambda_i$-tree $X_i$ for $i\in I$ and $\mathcal{D}$ is an ultrafilter in $I$ then the ultraproduct ${}^*G=\prod_{i\in I}G_i/\mathcal{D}$ has an essentially free ${}^*\theta$-affine action on the ${}^*\Lambda$-tree ${}^*X=\prod_{i\in I}X_i/\mathcal{D}$, where ${}^*\Lambda=\prod_{i\in I}\Lambda_i/\mathcal{D}$ and ${}^*\theta_{\langle g_i\rangle}\langle\lambda_i\rangle=\langle\theta_{g_i}^{(i)}\lambda_i\rangle$.
\end{proposition}

\pf
The action is as described in \cite[Theorem~39]{affine-paper}, where it is shown to be free.

For $\langle 1\rangle\neq \langle g_i\rangle\in{}^*G$, and for $i\in I$ let $\iota_i:A_g^{(i)}\to\Lambda_i$ be an isometric embedding.
We note first that $A_{\langle g_i\rangle}$ consists of $\langle x_i\rangle$ with $x_i\in A_{g_i}^{(i)}$ for almost all $i$: for $$\begin{array}{lll}\langle x_i\rangle\in A_{\langle g_i\rangle}&\Leftrightarrow {}^*d(\langle g_i\rangle^{-1}\langle x_i\rangle,\langle x_i\rangle)+{}^*d(\langle x_i\rangle,\langle g_i\rangle\langle x_i\rangle)={}^*d(\langle g_i\rangle^{-1}\langle x_i\rangle,\langle g_i\rangle\langle x_i\rangle)\\ &\Leftrightarrow \langle d_i(g_i^{-1}x_i,x_i)+d_i(x_i,g_ix_i)\rangle=\langle d_i(g_i^{-1}x_i,g_ix_i)\rangle\\ & \Leftrightarrow x_i\in A_{g_i}^{(i)}\ \mbox{for almost all } i.\end{array}$$

Define ${}^*\iota:A_{\langle g_i\rangle}\to{}^*\Lambda$ via $\langle x_i\rangle\mapsto\langle\iota_i(x_i)\rangle$. It is routine to show that ${}^*\iota$ is an isometric embedding of $A_{\langle g_i\rangle}$ in ${}^*\Lambda$. By Lemma~\ref{axis-embed} we now have ${}^*\iota(\langle g_i\rangle\langle x_i\rangle)={}^*\theta_{\langle g_i\rangle}{}^*\iota(\langle x_i\rangle)+\nu_{\langle g_i\rangle}$, giving $\langle\nu^{(i)}_{g_i}\rangle=\langle \iota_i(g_ix_i)-\theta^{(i)}_{g_i}(\iota_i(x_i))\rangle=\nu_{\langle g_i\rangle}$.
Now $(1-\theta_{\langle g_i\rangle})\langle \lambda_i\rangle = \langle (1-\theta_{g_i}^{(i)})(\lambda_i)\rangle
\ll\langle|\nu_{g_i}^{(i)}|\rangle=|\nu_{\langle g_i\rangle}|$. That is, $\nu_{\langle g_i\rangle}$ is ${}^*\theta_{\langle g_i\rangle}$-tame.\qed

\subsection{Regular embeddings of \oag s}

Given \oag s $\Lambda_0 $ and $\Lambda_1$, an $o$-embedding (that is, an order-preserving group embedding) $h$ of $\Lambda_0 $ in $\Lambda_1$, and an embedding $\bar{h}:\theta_\gamma\mapsto\eta_{\gamma}$ of $\Aut^+(\Lambda_0 )$ in $\Aut^+(\Lambda_1)$, we call $(h,\bar{h})$ an \emph{ample pair} for $\Lambda_0$ and $\Lambda_1$, and $h$ an \emph{ample embedding}
if \be
\item $h\cdot\theta_\gamma=\eta_{\gamma}\cdot h$ for all $\theta_\gamma\in\Aut^+(\Lambda_0 )$;
\item $[\im(1-\theta_g)]_{\Lambda_0 }\subset[\mu]_{\Lambda_0 }\Rightarrow[\im(1-\eta_g)]_{\Lambda_1}\subset[h(\mu)]_{\Lambda_1}$
\ee

The notation $\eta_g$ is intended to imply that $\eta_g=\bar{h}(\theta_g)$. In practice we will suppress $\bar{h}$, and leave the map to which it refers implicit in the rule $\theta_g\mapsto\eta_g$.

With regard to condition (b), note that it is automatic from (a) that $[\im(1-\theta_g)]\subset[\mu]$ implies $[\im(1-\eta_g)h]\subset[h(\mu)]$. If one also has $\im(1-\eta_g)]\subseteq[\im(1-\eta_g)h]$ then condition (b) follows; however in general the former condition need not be satisfied.

We will often suppress explicit mention of $\eta$
and refer to $h$ as an ample embedding.

A \emph{regular embedding} is an ample embedding for which
$\Aut^+(\Lambda_1)$ acts transitively on $h(\Lambda_0 ^{>0})$.
An \oag\ $\Lambda_0$ is \emph{regular} if $\Aut^+(\Lambda_0)$  acts transitively on $(\Lambda_0)^{>0}$.
A
\emph{strongly regular} embedding is an ample embedding
for which the codomain is a regular \oag.
Of course strongly regular embeddings are regular. It is
straightforward to show that the composition of ample, regular or
strongly regular embeddings is again ample, regular or strongly
regular, respectively.

\begin{lemma}\label{ample-embed}\be\item Let $\Lambda_\omega$ be \oag s where $\omega$ ranges through a linearly ordered set $\Omega$ and let $\Lambda=\omega(\prod_{\omega\in\Omega}\Lambda_\omega)$. Fix $\omega'\in\Omega$ and for a given $\theta_g^{\omega'}\in\Aut^+(\Lambda_{\omega'})$ define $\theta_g\in\Aut^+(\Lambda)$ via
$\theta_g(\lambda_\omega)_{\omega\in\Omega}=(\lambda'_\omega)_{\omega\in\Omega}$ where $\lambda'_\omega=\left\{\begin{array}{ll} \theta^{\omega'}_g\lambda_\omega'& \omega=\omega'\\ \lambda_\omega & \omega\neq\omega',\end{array}\right.$ and $h:\Lambda_{\omega'}\to\Lambda$ via $h(\lambda)=(\mu_\omega)_{\omega\in\Omega}$ where $\mu_\omega=\left\{\begin{array}{ll} \lambda & \omega=\omega'\\ 0 & \omega\neq\omega'.\end{array}\right.$

Then if $\bar{h}:\theta_g^{\omega'}\mapsto\theta_g$ then $(h,\bar{h})$ is an ample pair.

\item
Let $h:\Lambda_0 \to\Lambda_1$ be an ample embedding.
Then
if $G$ has an essentially free $\theta$-affine action on a $\Lambda_0 $-tree $X$, the induced action of $G$ on $X'=\Lambda_1\otimes_{\Lambda_0}  X$ is $\eta$-affine and essentially free.
\ee
\end{lemma}
\pf
(a) For $\lambda\in\Lambda_{\omega'}$, direct calculation shows that $\theta_g h(\lambda)=({\mu}_{\omega})_{\omega\in\Omega}=h\theta_g^{\omega'}(\lambda)$, where $$\mu_{\omega}=\left\{\begin{array}{cl}\theta_g^{\omega'}(\lambda) & \omega=\omega'\\ 0 & \omega\neq\omega'.\end{array}\right.$$ Next, suppose that $(1-\theta_g^{\omega'})(\lambda)\ll\mu$ for all $\lambda\in\Lambda_{\omega'}$. Then $(1-\theta_g)(\lambda_{\omega})_{\omega\in\Omega}=(\bar{\mu}_{\omega\in\Omega})$ where $$\bar{\mu}_\omega =\left\{\begin{array}{cl}(1-\theta_g^{\omega'})(\lambda_{\omega'}) & \omega=\omega'\\ 0 & \omega\neq\omega',\end{array}\right.$$ while the $\omega$th entry of $h(\mu)$ is equal to $\mu$ if $\omega=\omega'$ and 0 otherwise. The ampleness of $h$ follows.

(b)
We will use primes to refer to accessories of the action of $G$ on $X'$; thus, for example,
$A_g'$ denotes the axis of $g$ in $X'$ while $A_g$ is the axis of $g$ in $X$.

By \cite[Theorem 8(3)]{affine-paper}, there is a natural $\eta$-affine action of $G$ on $X'$ and a $G$-equivariant embedding $\phi:X\to X'$. Fix $g\neq 1$. Then $g$ is a hyperbolic automorphism of $X$, and for $x\in A_g\subseteq X$, using properties noted in this theorem, we have \begin{eqnarray*}d(g^{-1}x,gx)=d(g^{-1}x,x)+d(x,gx)
&\Rightarrow& h\ d(g^{-1}x,gx)=h\ d(g^{-1}x,x)+h\ d(x,gx)\\
&\Rightarrow& d'(\phi(g^{-1}x),\phi(gx))=d'(\phi(g^{-1}x),\phi(x))+d'(\phi(x),\phi(gx))\\
&\Rightarrow& d'(g^{-1}\phi(x),g\phi(x))=d'(g^{-1}\phi(x),\phi(x))+d'(\phi(x),g\phi(x))\\
&\Rightarrow& \phi(x)\in A_g'\subseteq X';
\end{eqnarray*}

Thus $\phi(A_g)\subseteq A_g'$. Now choose an isometric embedding
$\iota:A_g\hookrightarrow \Lambda_0 $. For distinct $x,y\in A_g$, we
have $\phi(x)\neq\phi(y)$, and we put $\iota'(\phi(x))=h\iota(x)$
and $\iota'(\phi(y))=h\iota(y)$.
Note that $|\iota'\phi(x)-\iota'\phi(y)|
=d'(\phi(x),\phi(y))$ and so,
by \cite[Lemma 2.3.1]{Chiswell-book} there is a
unique isometric extension of $\iota'$ to $A_g'$, and
$\iota'\phi=h\iota$.

Now by Lemma~\ref{ess-hyper-props}, $\iota(gx)=\theta_g\iota(x)+\nu_g$ for $x\in A_g$ and similarly for $A_g'$, giving
\begin{eqnarray*}\nu_g'&=&\iota'(g\phi(x))-\eta_g\iota'(\phi(x))\\ &=&\iota'(\phi(gx))-\eta_g\iota'(\phi(x))\\ &=&h\iota(gx)-\eta_g h\iota(x)\\ &=&h\iota(gx)-h\theta_g\iota(x)\\ &=&h\left(\iota(gx)-\theta_g\iota(x)\right)\\ &=&h(\nu_g)\end{eqnarray*}

Finally, essential freeness of the action on $X$ guarantees that for $g\neq 1$ we have $[\im(1-\theta_g)]\subset[\nu_g]$, while the ampleness of $h$ gives $[\im(1-\eta_g)]\subset[h(\nu_g)]=[\nu'_g]$;
thus $\nu_g'$ is $\eta_g$-tame.
\qed

\begin{proposition}\label{reg-embed}
\be\item Let $\Lambda_0 $ be an \oag\ of finite rank. Then $\Lambda_0 $ admits a strongly regular embedding in $\omega(\R^{\Z})$.
\item
Let $\Lambda_0 $ be an arbitrary \oag. There exists an \oag\ $\Lambda_1$ in which $\Lambda_0 $ admits a regular embedding.
\ee
\end{proposition}

\pf
(a) It is well-known that every \oag\ $\Lambda_0 $ of finite rank $n$ embeds in $\R^n$ (with the lexicographic order), whence it embeds in $\Lambda_1=\omega(\R^{\Z})$, via $h:(x_1,\ldots,x_n)\mapsto(x_m)_{m\in\Z}$, where $x_i=0$ if $i$ is not in the range between 1 and $n$. For $o$-automorphisms $\theta_g$ of $\R^n$, there is an extension to $\omega(\R^{\Z})$ defined by fixing the entries outside this range, which we will also denote by $\theta_g$.
By Lemma~\ref{ample-embed}(a), this embedding is ample.

Moreover, given positive $\lambda,\lambda^*\in\Lambda_1$, we can choose a shift $\sigma:(y_i)_{i\in\Z}\mapsto(y_{i-k})_{i\in\Z}$ so that $\sigma\lambda$ and $\lambda^*$ both have $i_0$ as the first position with a non-zero (and hence positive) entry. Applying the $o$-automorphism $\kappa$ that replaces the $i_0$th entry of $\sigma(\lambda)$ by a suitable scalar multiple, we can ensure that the first non-zero entries of $\kappa\cdot\sigma(\lambda)$ and $\lambda^*$ are equal. Finally, adding suitable multiples of the first non-zero entry to subsequent entries we obtain an $o$-automorphism $\xi$ of $\Lambda_1$ such that $\xi\cdot\kappa\cdot\sigma(\lambda)=\lambda^*$. Thus $\Aut^+(\Lambda_1)$ acts transitively on the positive elements of $\Lambda_1$.

(b) Let the finite rank subgroups of $\Lambda_0 $ be given by $\Lambda_i$ ($i\in I$). Let $\iota_i$ be a strongly regular embedding of $\Lambda_i$ in $\omega(\R^{\Z})$.

We can choose an ultrafilter $\mathcal{D}$ in $I$ such that $\Lambda_0 $ embeds in ${}^*\Lambda=\prod_{i\in I}\Lambda_i/\mathcal{D}$ via $\lambda\mapsto\langle \lambda_i\rangle$, say: one chooses $\mathcal{D}$ such that $\lambda=\lambda_i\in\Lambda_i$ for almost all $i$. Now form an embedding ${h}$ of $\Lambda_0$ in $\Lambda'=\prod_{i\in I}\omega(\R^{\Z})/\mathcal{D}$:

$${h}:\lambda\mapsto\langle \iota_i(\lambda_i)\rangle.$$ Defining $\eta_g:\langle \lambda_i\rangle\mapsto\langle \theta_g\lambda_i\rangle$ for $\theta_g\in\Aut^+(\Lambda_0)$, it is straightforward to check that $h$ is an ample embedding.
Moreover, for positive $\lambda=\langle\lambda_i\rangle,\lambda^*\langle\lambda_i^*\rangle\in\Lambda'$ and $i\in I$ using the regularity of $\iota_i$ we can obtain (for almost all $i$) a map $\theta_{g_i}\in\Aut^+\omega(\R^{\Z})$ such that $\theta_{g_i}(\iota_i\lambda_i)=\iota_i\lambda^*_i$. Thus $\theta_{\langle g_i\rangle}{h}(\lambda)=\theta_{\langle g_i\rangle}\langle \iota_i\lambda_i\rangle=\langle \theta_{g_i}(\iota_i\lambda_i) \rangle=\langle\iota_i\lambda^*_i\rangle={h}(\lambda^*)$; that is, $\Aut^+(\Lambda_1)$ acts transitively on $\im{h}$.
\h \qed

\subsection{Some combinations of $\ATF$ groups}

Suppose now that a signature (S1)-(S5) is given as in Theorem~\ref{Bass3.8'} and consider the following condition.\\

(IE) $\theta_{\alpha_e(s)}^{\partial_0e}=1$ for all $s\in\mathcal{G}(e)$, $e\in E(Y^*)$. \\

We will say that \emph{edge groups are isometric} in this case (note that this shorthand does not imply that $\alpha_e(\mathcal{G}(e))$ is contained in the kernel of $\beta$). Then assuming (C1) we have end homomorphisms $\tau_e:\alpha_e\mathcal{G}(e)\to\Lambda_0$ such that $\ell(\alpha_e(g))=|\tau_e\alpha_e(g)|$ ($g\in\mathcal{G}(e)$). See \cite[1.6]{Bass} for the definition and basic properties of end homomorphisms. Note that if in the context of Theorem~\ref{Bass3.5'} the action of $\Gamma$ on a $\Z\times\Lambda_0$-tree is free and (IE) is satisfied, then the action of each $\alpha_e\mathcal{G}(e)$ on $X(\partial_0e)$ is free and isometric and fixes the end $\epsilon_e$. Thus there is a line $(\epsilon_e,\bar{\epsilon}_{{e}})$ of $X(\partial_0e)$ stabilised by $\alpha_e\mathcal{G}(e)$. (Recall that a line in a $\Lambda_0$-tree $X_0$ is a maximal linear subset of $X_0$.)

In the presence of assumption (IE), conditions (C3) and (C4)
are respectively equivalent to

\be\item[($\Cthree_{IE}$)] $\theta_{g_e}=1$ if $g_e=1$;
\item[($\Cfour_{IE}$)]  $\theta_{g_e}\mu_{\alpha_e(s)}=\mu_{\alpha_{|e|}(s)}$ ($s\in\mathcal{G}(e)$, $e\in E(Y^*)$), and $\mu^{x^*}:\mathcal{G}(x^*)\to\Lambda_0$ is a homomorphism ($x^*\in Y^*$). Moreover, $\mu_{g_e}=0$ if $g_e=1$;
\ee
Furthermore, assuming (IE), we have

\begin{eqnarray*} \theta_{g_e}\left[\theta_{\alpha_e(s)}\delta_e-\delta_e\cdot \alpha_e(s)\right]+\theta_{g_{\bar{e}}}\left[\theta_{\alpha_{\bar{e}}(s)}\delta_{\bar{e}}-\delta_{\bar{e}}\cdot \alpha_{\bar{e}}(s)\right]&=& \theta_{g_e}\left[\delta_e-\delta_e\cdot \alpha_{e}(s)\right]+\theta_{g_{\bar{e}}}\left[\delta_{\bar{e}}-\delta_{\bar{e}}\cdot \alpha_{\bar{e}}(s)\right]\\ &=& -\left[\theta_{g_e}\cdot\tau_e\cdot\alpha_e(s)+\theta_{g_{\bar{e}}}\cdot\tau_{\bar{e}}\cdot\alpha_{\bar{e}}(s)\right]\end{eqnarray*}
and so the compatibility condition (C5) is equivalent to the tidier identity\\

\be\item[($\Cfive_{IE}$)] $\theta_{g_e}\cdot\tau_e\cdot\alpha_e(s)+\theta_{g_{\bar{e}}}\cdot\tau_{\bar{e}}\cdot\alpha_{\bar{e}}(s)=-\mu_{\alpha_{|e|}(s)},\h s\in\mathcal{G}(e)$.\ee

We will specialise further in several of the following corollaries to the case where the edge groups are cyclic. That is, we will assume

\be\item[(CIE)] For $e\in E(Y^*)$ we have $\mathcal{G}(e)$ cyclic and $\theta_{\alpha_e(s)}^{\partial_0^{Y^*}e}=1$ for $s\in\mathcal{G}(e)$.\ee

Now consider the following weakened version of (C1).

\be\item[($\Cone_{\wedge}$)] $\alpha_e\mathcal{G}(e)=(\mathcal{G}(\partial_0e))_{\epsilon_e}$ ($e\in E(Y^*)$).\ee

\begin{corollary}\label{cor-isom-edge-gps}
Suppose that data as in (S1)-(S5) are given (see Theorem~\ref{Bass3.8'}).
\be\item Assume that
(IE), and the compatibility conditions (C1), (C2), ($\mathit{C3}_{IE}$), ($\mathit{C4}_{IE}$) and ($\mathit{C5}_{IE}$) are satisfied. Then $\pi_1(\mathcal{G})$ admits an affine action as in (A1)-(A4). If in addition the actions of the vertex groups are free, so is the resulting action of $\pi_1(\mathcal{G})$.
\item
Assume that the compatibility conditions ($\mathit{C1}_{\wedge}$) and (C2)-(C5) are satisfied. Then $\pi_1(\mathcal{G})$ admits an affine action as in (A1)-(A4). If the vertex group actions are essentially free, then the action of $\pi_1(\mathcal{G})$ is free.

\item
If the given actions of the vertex groups are essentially free and condition (IE) is satisfied, and if ($\mathit{C1}_{\wedge}$), (C2), ($\mathit{C3}_{IE}$), ($\mathit{C4}_{IE}$), ($\mathit{C5}_{IE}$) then $\pi_1(\mathcal{G})$ admits a free affine action as in (A1)-(A4).\ee
\end{corollary}

\pf
Part (a) summarises the observations above, while the assertions of
parts (b) and (c) are easy consequences of Lemma~\ref{fulfilment} and Lemma~\ref{ess-hyper-props}(d).
\qed

\begin{lemma}\label{rescale-metric}
Let $\Gamma_0$ be a group equipped with a $\theta$-affine action on a $\Lambda_0$-tree $(X_0,d_0)$. Let $\epsilon$ be an end of $X_0$, and $\delta:X_0\to\Lambda_0$ an end map towards $\epsilon$. Let $\eta\in\Aut^+(\Lambda_0)$ and $d_1=\eta\cdot d_0$. Then $(X_0,d_1)$ is a $\Lambda_0$-tree, the action of $\Gamma_0$ on $(X_0,d_1)$ is $(\eta\cdot\theta\cdot\eta^{-1})$-affine, and the end maps with respect to the metrics $d_0$ and $d_1$ are in natural bijective correspondence via $\delta_0\mapsto\eta\cdot\delta_0$.
If $\epsilon$ is an end of full $\Lambda_0$-type with respect to $d_0$, it is of full $\Lambda_0$-type with respect to $d_1$.
\end{lemma}

\pf
Writing $[x,y]_i$ for segments in $X_i$, one notes first that
$d_0(x,z)=d_0(x,y)+d_0(y,z)$ if and only if $\eta d_0(x,z)=\eta d_0(x,y)+\eta d_0(y,z)$; it follows that
$[x,z]_0=[x,z]_1$.
It follows that $(X_0,d_1)$ is indeed a $\Lambda_0$-tree with the same segments, $X_0$-rays and ends as $(X_0,d_0)$. Moreover if $\epsilon$ is of full $\Lambda_0$-type with respect to $d_0$ then $\{d_0(x,y):y\in[x,\epsilon)\}=[0,\infty)_{\Lambda_0}$ giving $\{d_1(x,y):y\in[x,\epsilon)\}=\eta[0,\infty)_{\Lambda_0}=[0,\infty)_{\Lambda_0}$, whence $\epsilon$ is of full $\Lambda_0$-type with respect to $d_1$; the converse is obvious.

Now if $\delta_0$ is an end map towards $\epsilon$ with respect to $d_0$ then $\delta_0(y)-\delta_0(x)=d_0(x,y)$ for $y\in[x,\epsilon)$, whence $\eta\cdot\delta_0(y)-\eta\cdot\delta_0(x)=\eta\cdot d_0(x,y)$; that is, $\delta_1=\eta\cdot\delta_0$ is an end map with respect to $d_1$. Reversing this argument gives the required assertion concerning end maps.

That the action on $(X_0,d_1)$ is $(\eta\cdot\theta\cdot\eta^{-1})$-affine is a straightforward calculation.
\qed

An easy calculation also shows that $\theta_{g_e}\cdot\tau_e\cdot\alpha_e(st)=\theta_{g_e}\cdot\tau_e\cdot\alpha_e(s)+\theta_{g_e}\cdot\tau_e\cdot\alpha_e(t)$. Thus if $\mathcal{G}(e)=\langle s_e\rangle$ is cyclic, and the equation in ($\Cfive_{IE}$) holds for $s=s_e$, it follows that ($\Cfive_{IE}$) holds for all $s\in\mathcal{G}(e)$.

\begin{lemma}\label{ess-mu-g}
Let $g$ be a hyperbolic $\beta_g$-affine automorphism of a $\Lambda$-tree $X$ (where $\Lambda=\Z\times\Lambda_0$), and suppose that $g$ stabilises some ball $X_0$ of radius $0\times\Lambda_0$. Suppose that $\beta_g:(1,\lambda_0)\mapsto(1,\theta_g\lambda_0+\mu_g)$.

Let $\bar{\iota}:\bar{A_g}\to\Lambda_0$ be an isometric embedding, where $\bar{A_g}$ is the axis of $g$ with respect to the action of $g$ on $X_0$ (thus $\bar{A_g}$ is equal to $A_g\cap X_0$, and is viewed as a $\Lambda_0$-tree). There exists an isometric embedding $\iota:A_g\to\Lambda$ such that $\iota(x)=(0,\bar{\iota}(x))$ for $x\in\bar{A}_g$. If $\iota(gx)=\beta_g \iota(x)+\nu_g$ and $\bar{\iota}(gx)=\theta_g \bar{\iota}(x)+\bar{\nu}_g$ then $\nu_g=(0,\bar{\nu}_g)$.

If $g$ is essentially hyperbolic as an automorphism of $X_0$, then $g$ is an essentially hyperbolic automorphism of $X$ if and only if $\mu_g\ll|\bar{\nu}_g|$.
\end{lemma}

\pf
Put $\hat{\iota}(x)=(0,\bar{\iota}(x))$; this defines a $\Lambda$-isometric embedding of $\bar{A}_g$ in $\Lambda$. Since $g$ is hyperbolic, $\bar{A}_g$ contains more than one point, and we may apply \cite[Lemma 2.3.1]{Chiswell-book} to extend $\hat{\iota}$ to an embedding $\iota$ of $A_g$ in $\Lambda$.

By Lemma~\ref{axis-embed} there exists $\nu_g$ such that for $x\in \bar{A}_g$
we have

\begin{eqnarray*}\iota(gx)&=&\beta_g\iota(x)+\nu_g\\ &=& (0,\theta_g\bar{\iota}(x))+\nu_g.\end{eqnarray*}

Note that $gx\in \bar{A}_g$ and therefore \begin{eqnarray*}\iota(gx)&=&(0,\bar{\iota}(gx))\\ &=& (0,\theta_g\bar{\iota}(x)
+\bar{\nu}_g),\end{eqnarray*}
\\
so comparing these two expressions, we get $\nu_g=(0,\bar{\nu}_g)$.

Now suppose that $(1-\beta_g)(k,\lambda_0)\ll|\nu_g|$ for all $(k,\lambda_0)\in\Lambda$. Then $(0,(1-\theta_g)(\lambda_0))\ll |(0,\bar{\nu}_g)|$ and  $(0,(1-\theta_g)(\lambda_0)-k\mu_g)\ll|(0,\bar{\nu}_g)|$, whence $k\mu_g\ll|\bar{\nu}_g|$ for all $k$; in particular, $\mu_g\ll|\bar{\nu}_g|$.

Conversely, suppose that $\mu_g\ll|\bar{\nu}_g|$. Now the essential hyperbolicity of $g$ (with respect to $X_0$) gives $(1-\theta_g)(\lambda_0)\ll|\bar{\nu}_g|$ for all $\lambda_0\in\Lambda_0$, whence $(1-\beta_g)(k,\lambda_0)=(0,(1-\theta_g)(\lambda_0)-k\mu_g)\ll|(0,\bar{\nu}_g)|=|\nu_g|$.
\qed

Consider now the following condition, a modification of (C2).\\

($\Ctwo'$) If $\partial_0e_i=x^*$ ($i=1,2,3$) and $\mathcal{G}(e_1)\neq
1$, then the $\alpha_{e_i}\mathcal{G}(e_i)$ are not all conjugate in
$\mathcal{G}(x^*)$.\\

\begin{lemma}\label{edge-ends}
Suppose that a graph of groups $(\mathcal{G},Y^*,T^*)$ is given as in (S1), with free $\theta^{x^*}$-affine actions of $\mathcal{G}(x^*)$ on $\Lambda_{x^*}$-trees, and that (CIE) and (C2') are satisfied.
Suppose that $\alpha_e\mathcal{G}(e)$ coincides with the stabiliser in $\mathcal{G}(\partial_0e)$ of some end of $X(\partial_0e)$ and that $\alpha_e(s)$ is not conjugate in $\pi_1(\mathcal{G})$ to $\alpha_e(s^{-1})$ whenever $e\in E(Y^*)$ and $1\neq s\in\mathcal{G}(e)$. \be\item There exists an orientation $E^+\subseteq E(Y^*)$ and, for $e\in E(Y^*)$ with $\mathcal{G}(e)\neq 1$, there are generators $s_e$ of $\mathcal{G}(e)=\mathcal{G}(\bar{e})$ such that $\alpha_e(s_e)$ translates towards $\epsilon_e$ and $\alpha_{\bar{e}}(s_e^{-1})$ translates towards $\epsilon_{\bar{e}}$ if $e\in E^+$.
\item For $e\neq f$ with $\partial_0e=\partial_0f=x^*$ the ends $\epsilon_e$ and $\epsilon_f$ are in distinct $\mathcal{G}(x^*)$-orbits.
\ee
\end{lemma}

\pf
Consider an edge $e$ for which the edge group $\mathcal{G}(e)$ is non-trivial, and define an oriented subgraph $Y_e^*$ as follows. Declare $e=e_0$ and $\bar{e}=\bar{e}_0$ to be edges of $Y_e^*$ with $e$ as the positively oriented edge, and inductively declare $e_{k+1}$ and $\bar{e}_{k+1}$ to be edges of $Y_e^*$ (with $e_{k+1}$ as the positively oriented edge) if $\partial_0e_{k+1}=\partial_0\bar{e}_{k}$ and  $\alpha_{e_{k+1}}\mathcal{G}(e_{k+1})$ is conjugate (in $\mathcal{G}(\partial_0\bar{e}_k)$) to $\alpha_{\bar{e}_k}\mathcal{G}(e_k)$ for some $e_k\in Y_e^*$. It is easy to see that $Y_e^*$ is connected and that the embedded images $\alpha_{e'}\mathcal{G}(e')$ of all edge groups in $\pi_1(\mathcal{G})$ are conjugate and non-trivial for $e'\in E(Y_e^*)$. We claim that the degree of each vertex of $Y_e^*$ is at most 2. For otherwise,
if $x^*$ has degree greater than 2, there exist $f_0,f_1,f_2\in E(Y^*)$ such that $\partial_0f_i=x^*$ for each $i$.
Then the subgroups $\alpha_{f_i}\mathcal{G}(f_i)$ are all conjugate in $\mathcal{G}(x^*)$ and non-trivial (since $\mathcal{G}(e)\neq 1$), which violates condition (C2').

It is also straightforward to see that $Y_{e_{k+1}}^*=Y_{e_k}^*$ for $e_{k}\in Y_e^*$. It follows that the set of edges with non-trivial edge group is expressible as a disjoint union of sets of the form $E(Y_{f}^*)$.

For each $Y_e^*$ take a positively oriented edge $e_0\in E(Y_e^*)$ (possibly the same $e_0$ as above) and choose a generator $s_{e_0}=s_{\bar{e}_0}\in \mathcal{G}(e_0)$. Then there are exactly two ends of $X(\partial_0e_0)$ fixed by $\alpha_{e_0}(s_{e_0})$: let $\epsilon_{e_0}$ be the end
for which $\alpha_{e_0}(s_{e_0})$ translates towards $\epsilon_{e_0}$, and denote the other end by $\bar{\epsilon}_{e_0}$. We let $\epsilon_{\bar{e}_0}$ be the end of the axis of $\alpha_{\bar{e}_0}(s_{e_0})$ (in $X(\partial_0\bar{e}_0)$) for which $\alpha_{\bar{e}_0}(s_{e_0})$ translates away from $\epsilon_{\bar{e}_0}$, and denote the other end of this axis by $\bar{\epsilon}_{\bar{e}_0}$. Inductively, if $\partial_0\bar{e}_k$ is incident to another edge $e_{k+1}$ of $Y_e^*$, so that $\partial_0e_{k+1}=\partial_0\bar{e}_k$, $e_{k+1}\neq\bar{e}_k$, $\partial_0\bar{e}_{k+1}$ is not incident to any edge of $Y_e^*$, and
$u\alpha_{\bar{e}_k}\mathcal{G}(e_k)u^{-1}=\alpha_{e_{k+1}}\mathcal{G}(e_{k+1})$, then let the generator $s_{e_{k+1}}$ be chosen such that $\alpha_{e_{k+1}}(s_{e_{k+1}})=u\alpha_{\bar{e}_k}(s_{\bar{e}_k})u^{-1}$. Since $\alpha_{e_{k+1}}(s_{e_{k+1}})$ is not conjugate to its inverse, this choice of $s_{e_{k+1}}$ is independent of $u$. Moreover $\alpha_{\bar{e}_k}(s_{e_k})$ translates away from $\epsilon_{\bar{e}_k}$, and hence $\alpha_{e_{k+1}}(s_{e_{k+1}})=u\alpha_{\bar{e}_k}(s_{e_k})u^{-1}$ translates away from $u\epsilon_{\bar{e}_k}$.
We let $\epsilon_{e_{k+1}}$ be the other end of the axis of $\alpha_{e_{k+1}}(s_{e_{k+1}})$ (in $X(\partial_0e_{k+1})$). Thus $\alpha_{e_{k+1}}(s_{e_{k+1}})$ translates towards $\epsilon_{e_{k+1}}$. Since $\bar{e}_k$ and $e_{k+1}$ are the only edges of $Y_e^*$ incident to $\partial_0e_{k+1}$, our choice of $\epsilon_{e_{k+1}}$ ensures that condition (b) is satisfied.

For edges $e_{-1}\in E(Y_e^*)$ with $\partial_0e_0=\partial_0\bar{e}_{-1}$ incident to $e_{-1}$ and $e_{-1}\notin\{e_k,\bar{e}_k:k\geq 1\}$ we can define $s_{e_{-1}}$ and ends $\epsilon_{e_{-1}}$, and $e_{-k}$ and $s_{e_{-k}}$ in a similar fashion.

We thus obtain ends $\epsilon_{e_k}$ for each $e_k\in E(Y_e^*)$ in a maximal subtree of $Y_e^*$
and the required conditions (a) and (b) are satisfied for edges of this subtree. If $Y_e^*$
is not a tree there is exactly one further edge $e_m\in Y_e^*$, with $\partial_0e_m=\partial_0\bar{e}_{m-1}$ and $\partial_0\bar{e}_m=\partial_0e_0$. Now $\alpha_{\bar{e}_m}\mathcal{G}(e_m)$ is conjugate to $\alpha_{e_0}\mathcal{G}(e_0)$. In fact, the Bass-Serre relations for $\pi_1(\mathcal{G})$ ensure that $\alpha_{e_k}(s_k)$ is conjugate to $\alpha_{\bar{e}_k}(s_k)$ for all $k$, and the definition of $Y_e^*$ ensures that $\alpha_{e_k}(s_k)$ is conjugate to $\alpha_{e_{k+1}}(s_{k+1})$ for $0\leq k\leq m-1$. Thus $\alpha_{e_0}(s_{e_0})$ is conjugate to $\alpha_{\bar{e}_m}(s_{e_m})$. Our assumption that no non-trivial element of an embedded edge group can be conjugate in $\Gamma$ to its inverse ensures that the ends $\epsilon_{e_m}$ and $\epsilon_{\bar{e}_m}$ can be chosen so that $\alpha_{e_m}(s_{e_m})$ translates towards $\epsilon_{e_m}$ and $\alpha_{\bar{e}_m}(s_{e_m})$ translates away from $\epsilon_{\bar{e}_m}$.

Repeating this procedure for the other graphs $Y_{e'}^*$ results in an assignment $e\mapsto\epsilon_e$
with the required properties (a) and (b) for all edges with non-trivial edge group.
For trivial $\mathcal{G}(e)$ we
choose $\epsilon_e$ arbitrarily among the ends that have trivial
stabiliser in such a way that
condition (b) is satisfied. (That
such an assignment of ends to edges is possible in this case is essentially shown
in \cite[4.9]{Bass}.)

Finally suppose that $\epsilon_e=\gamma\epsilon_f$ for edges $e$ and $f$ with $\partial_0e=\partial_0f=x^*$ and non-trivial $\gamma\in\mathcal{G}(x^*)$. Then $\mathcal{G}(e),\mathcal{G}(f)\neq 1$, and $\alpha_e(s_e^{\chi_e})$ and $\gamma\alpha_f(s_f^{\chi_f})\gamma^{-1}$ both translate towards $\gamma\epsilon_f$ where $\chi_{e'}=1$ if $e'\in E^+$ and $\chi_{e'}=-1$ if $\bar{e'}\in E^+$. Since both of these elements generate end stabilisers of $\mathcal{G}(x^*)$, we must have $\gamma\mathcal{G}(f)\gamma^{-1}=\alpha_e\mathcal{G}(e)$. Thus $Y_e^*=Y_f^*$. But we have already established that if $f\in Y_e^*$ then $\epsilon_e$ and $\epsilon_f$ lie in distinct $\mathcal{G}(x^*)$ orbits for $e\neq f$; thus $e=f$.
\qed

We can now establish the main result of this section which concerns combinations of groups that admit essentially free actions. The proof exploits the degrees of freedom afforded by the $\theta_{g_e}$ to ensure that under suitable hypotheses the condition (C5) is satisfied with $\mu_g=0$ for all $g$.

\begin{theorem}\label{relax-C5-thm}
Let a graph of groups $(\mathcal{G},Y^*,T^*)$ be given as in (S1), and for $x^*\in Y^*$ let a $\Lambda_{x^*}$-tree $X(x^*)$ be given on which $\mathcal{G}(x^*)$ has an essentially free $\bar{\theta}^{x^*}$-affine action (where $\bar{\theta}^{x^*}:\mathcal{G}(x^*)\to\Aut^+(\Lambda_{x^*})$ is a homomorphism). For $e\in E(Y^*)$ suppose that $\alpha_e\mathcal{G}(e)$ coincides with the stabiliser in $\mathcal{G}(\partial_0e)$ of some end of $X(\partial_0e)$,
and that $\alpha_e(g)$ is not conjugate in $\pi_1(\mathcal{G})$  to $\alpha_e(g^{-1})$ for $1\neq g\in\mathcal{G}(e)$.

Assume further that conditions (CIE) and (C2') are satisfied.
Then $\Gamma=\pi_1(\mathcal{G})$ admits an essentially free affine action on a $\Lambda$-tree.

If $\Lambda_{x^*}=\Lambda_0$ for all $x^*\in Y^*$ and $\Lambda_0$ is regular we can take $\Lambda=\Z\times\Lambda_0$.
\end{theorem}

\pf
We first arrange things so that the \oag s $\Lambda_{x^*}$ match up. We can do this by taking an arbitrary linear order on $Y^*$ and embedding each $\Lambda_{x_0^*}$ in $\Lambda_1=\omega(\prod_{x^*\in Y^*}\Lambda_{x^*})$ as in Lemma~\ref{ample-embed}(a). Furthermore, using Proposition~\ref{reg-embed}(b), we can find a regular embedding $h$ of $\Lambda_1$ in an \oag\ $\Lambda_0$. Applying the base change functor to each of these embeddings, we therefore lose no generality in assuming that $\Lambda_{x^*}=\Lambda_0$ for all $x^*\in Y^*$. Moreover the actions of the vertex groups on the resulting $\Lambda_0$-trees are also essentially free by Lemma~\ref{ample-embed}. Of course if $\Lambda_{x^*}=\Lambda_0$ for all $x^*$ where $\Lambda_0$ is regular, these adjustments to $\Lambda_{x^*}$ are unnecessary.

We now use Lemma~\ref{fulfilment} to extend the actions to the respective $\Lambda_0$-fulfilments, so all ends are of full $\Lambda_0$-type and the resulting actions are still free, by Lemma~\ref{ess-hyper-props}(d).

We now assign an end $\epsilon_e$ to each $e\in E(Y^*)$ as in Lemma~\ref{edge-ends}. Then the quantities $\tau_e\cdot\alpha_e(s)=\delta_e\cdot\alpha_e(s)-\delta_e$ and $\tau_{\bar{e}}\cdot\alpha_{\bar{e}}(s)=\delta_{\bar{e}}\cdot\alpha_{\bar{e}}(s)-\delta_{\bar{e}}$ have opposite sign for non-trivial $s\in\mathcal{G}(e)$. That both are constants is shown in \cite[1.4(b)]{Bass}.

Take a vertex $x_0^*$ of $Y^*$, and for each edge $e\in E(T^*)$ with $\partial_0^{Y^*}\bar{e}=x_0^*$, choose $\eta_e\in\Aut^+(\Lambda_0)$ such that $$\eta_e\cdot\tau_e\cdot\alpha_e(s_e)+\tau_{\bar{e}}\cdot\alpha_{\bar{e}}(s_e)=0.$$
(Since $\Aut^+(\Lambda_0)$ acts transitively on the embedded image of $\Lambda_1$ and $\tau_e\cdot\alpha_e(s_e)\in\Lambda_1$ for all $e$, the existence of such an $\eta_e$ follows from the regularity of the embedding of $\Lambda_1$ in $\Lambda_0$ and our choice of $s_e$ which ensures that $\tau_e\alpha_e(s_e)$ and $\tau_{\bar{e}}\alpha_{\bar{e}}(s_{\bar{e}})$ have opposite signs.)

Since $\mathcal{G}(e)$ is cyclic, it follows that the same identity holds with $s_e$ replaced by any $t\in\mathcal{G}(e)$.

Put $x_1^*=\partial_0^{Y^*}e$ and replace the given $\Lambda_0$-metric $d_{x_1^*}$ on $X(x_1^*)$ by $\eta_e\cdot d_{x_1^*}$. Performing this replacement for all edges $e$ incident to $x_0^*$ and applying Lemma~\ref{rescale-metric}, one now has condition ($\Cfive_{IE}$) satisfied for each $e\in E(T^*)$ incident to $x_0^*$ at the expense of replacing the homomorphisms $\bar{\theta}^{x^*}$ by $\theta^{x^*}=\eta_e\cdot\bar{\theta}^{x^*}\cdot\eta_e^{-1}$.
Since we have restricted to edges of a subtree of $Y^*$, these adjustments to the metrics on the $\Lambda_0$-trees do not come into conflict. Inductively, we can continue this process to adjust the metric on each $X(x_k^*)$ where the length of the reduced path in $T^*$ joining $x_k^*$ to $x_0^*$ is $k$. We then have ($\Cfive_{IE}$) satisfied for all $e\in E(T^*)$.

For each edge $e\in E(Y^*)$ with $\partial_0^{X^*}e\notin T^*$ we
now choose $\theta_{g_e}\in\Aut^+(\Lambda_0)$ such that
$$\theta_{g_e}\cdot\tau_e\cdot\alpha_e(s_e)+\tau_{\bar{e}}\cdot\alpha_{\bar{e}}(s_e)=0.$$
For $e$ with $\partial_0^{X^*}e\in T^*$ we put $\theta_{g_e}=1$.
Setting $\mu_{g_e}=\mu_{s}^{x^*}=0$ for all $e$ and all $s\in\mathcal{G}(x^*)$, and taking an arbitrary end map $\delta_e$ towards $\epsilon_e$ for $e\in E(Y^*)$, we now have all the data described by (S1)-(S5) of Theorem~\ref{Bass3.8'}, and we have ensured that these data satisfy (C1), (C2), ($\Cthree_{IE}$), ($\Cfour_{IE}$) and ($\Cfive_{IE}$).

The existence of a free $\beta$-affine action of $\Gamma$ now
follows from Corollary~\ref{cor-isom-edge-gps}. To show that the action is essentially free, consider $g\in\Gamma$. If $g$ is not in a conjugate of a vertex group then $g$ stabilises no ball of radius $0\times\Lambda_0$, whence $\nu_g=(m,\kappa)$ where $m\neq 0$. Thus $(1-\beta_g)(m,\lambda_0)\ll|\nu_g|$ in this case.

Otherwise some conjugate $g_0$ of $g$ belongs to a vertex group. That such a $g_0$ --- and hence such a $g$ --- is essentially hyperbolic follows from Lemma~\ref{ess-mu-g}.
\qed

The following definition appears in \cite[\S1.4]{affine-paper}. Let
$\mathcal{P}$ be a class of groups (with $\{1\}\in\mathcal{P}$). A group $G$ is said to be
$\ATF[\mathcal{P}]$ if $G$ admits a free affine action on a
$\Lambda_0 $-tree (for some $\Lambda_0 $) such that all line stabilisers
are $\mathcal{P}$~subgroups, and every $\mathcal{P}$ subgroup of $G$
stabilises a line. In this context we will call a subgroup $H$ of an
$\ATF[\mathcal{P}]$ group $G$ an \emph{end subgroup} if either $H$ is trivial and $G$ is not $\mathcal{P}$, or $H$
is maximal $\mathcal{P}$ in $G$ -- we include here the possibility
that $H=G$. This
amounts to the assertion that $H=G_{\epsilon}$, the stabiliser of
some end $\epsilon$ (with respect to a minimal action). Bass defines end subgroups in the context of
isometric actions on $\Lambda_0 $-trees, where all end stabilisers are
automatically abelian. We will consider here the case where
$\mathcal{P}=\mathrm{sol}$, the class of soluble groups. This contains all ITF groups as well as $\ATF(\Lambda_0)$ groups for which $\Lambda_0$ has finite rank. Note that end stabilisers coincide with line stabilisers in this case; see \cite[Proposition 27]{affine-paper}.

\begin{notation}\label{notation} We will further encumber the notation above with a
superscript $o$ if in addition the action preserves the orientation of each line. Similarly we may add a superscript
$e$ if in addition the action can be taken to be
essentially free. For example, $\ATF^o(\Lambda_0)$ refers to groups that admit a free (but not necessarily essentially free) affine action on a $\Lambda_0$-tree where no line has its ends interchanged by any group element; $\ATF^e[\mathrm{sol}]$ refers to groups $G$
that admit an essentially free affine action with end stabilisers
that are either trivial or maximal soluble in $G$.
\end{notation}

\begin{lemma}\label{line-stabs-BCF}
Let $(X_0,d_0)$ be a $\Lambda_0$-tree, and $h:\Lambda_0\to\Lambda_1$ an $o$-embedding. Let $(X_1,d_1)$ and embedding $\phi:X_0\to X_1$ be the result of applying the base change functor to $(X_0,d_0)$ and $h$. Suppose that $\beta:G\to\Aut^+(\Lambda_0)$ and $\eta:G\to\Aut^+(\Lambda_1)$ are homomorphisms where $h\beta_g=\eta_g h$, and $G$ has a $\beta$-affine action on $X_0$. If $L_0$ is a line in $X_0$, then $\phi(L_0)$ spans a line $L_1$ in $X_1$. If $H$ is the stabiliser in $G$ of $L_0$, then
$H$ is the stabiliser in $G$ of $L_1$.
\end{lemma}

\pf
We will use the basic properties of the base change functor as in \cite[Theorem 8(3)]{affine-paper}.
Suppose that $L_0$ is a line in $X_0$. If $x,y,z\in L_0$ with $y\in[x,z]_0$ then $d_0(x,z)=d_0(x,y)+d_0(y,z)$. Applying $h$ to both sides of this equation and observing that $d_1(\phi(u),\phi(v))=h d_0(u,v)$ we see that $\phi(y)\in[\phi(x),\phi(z)]_1$. Thus $\phi(L_0)$ is a linear subset of $X_1$. To see that it spans a maximal linear subset of $X_1$, suppose that $w\in X_1$ satisfies $[\phi(x),w]_1\supset[\phi(x),\phi(u)]_1$ for all $u\in L_0$.
Since $X_1$ is spanned by $\phi(X_0)$, we can take $v\in X_0$ with $w\in[\phi(x),\phi(v)]_1$. Now $[\phi(x),w]_1\subseteq[\phi(x),\phi(v)]_1$.
Letting $\epsilon$ and $\epsilon'$ denote the ends of $L_0$ in $X_0$ and replacing $v$ by $\Y(v,\epsilon,\epsilon')$, we have $v\in L_0$, whence $\phi[x,v]_0\subset[\phi(x),w]_1$, since $[\phi(x),w]_1$ is a full subset of $X_1$ containing $\phi(x)$ and $\phi(v)$. But this contradicts $[\phi(x),w]_1\subseteq[\phi(x),\phi(v)]_1$. Thus $\phi(L_0)$ spans a line in $X_1$, which we will denote by $L_1$.

Suppose now that $H$ is the stabiliser of $L_0$. Since $\phi$ is an equivariant embedding
$H$ must stabilise $\phi(L_0)$, and thus $L_1$. Conversely, suppose that $gL_1=L_1$, and that $x\in L_0$, so that $\phi(x)$ and $\phi(gx)=g\phi(x)$ are in $L_1$. Choose $u,v\in L_0$ such that $\phi(x),\phi(gx)\in[\phi(u),\phi(v)]_1$. Then since $h$ is injective we have $x,gx\in[u,v]_0$. It follows that $g\in H$. Therefore $H$ is the stabiliser of $L_1$, as claimed.
\h\qed

\begin{lemma}\label{sol-line-stabs}
Suppose that $\Gamma$ has a free $\beta$-affine action on the $\Lambda$-tree $X$ (where $\Lambda=\Z\times\Lambda_0$), and that the actions of $\Gamma_{x^*}$ on the balls $X(x^*)$ of radius $0\times\Lambda_0$ are free with soluble line stabilisers. Then the action of $\Gamma$ has soluble line stabilisers.
\end{lemma}

\pf
Suppose that $H\neq 1$ is the stabiliser of the line $L\subseteq X$, and
note that $N=H\cap \Gamma_{x^*}$ stabilises the line $L\cap X(x^*)$ in $X(x^*)$, viewed as a $0\times\Lambda_0$-tree, for some $x^*\in X^*$. Our assumption on the vertex groups implies that $N$ is soluble. Moreover $N$ is the kernel of the action of $H$ on the linear subtree of $X^*$ stabilised by $H$; thus $N$ is normal in $H$ with cyclic quotient $H/N$ whence $H$ is soluble.
\qed

\begin{corollary}\label{C2'-cor}
Let a graph of groups $(\mathcal{G},Y^*,T^*)$ be given as in (S1), where each
vertex group $\mathcal{G}(x^*)$ has an essentially free affine action on a $\Lambda_0$-tree with soluble line stabilisers.
For $e\in E(Y^*)$ suppose that
$\alpha_e\mathcal{G}_e$ is an end subgroup
of $\mathcal{G}(\partial_0e)$.

Assume that
conditions (CIE) and (C2') are satisfied, and that no non-trivial element of $\Gamma$ is conjugate to its inverse.

Then $\Gamma$ is $\ATF^e[\mathrm{sol}]$.

\end{corollary}

\pf
Let essentially free $\bar{\theta}^{x^*}$-affine actions of the vertex groups $\mathcal{G}(x^*)$ on $\Lambda_{x^*}$-trees be given where line stabilisers are soluble. It follows from Lemma~\ref{line-stabs-BCF} that if $H$ stabilises the line $L$ with respect to the affine action of a group on a $\Lambda_1$-tree, then the same is true of the induced action of $H$ on the $\Lambda_2$-tree obtained by the Base Change Functor applied to an embedding $h:\Lambda_1\to\Lambda_2$. We can thus adjust the original affine actions as in the proof of Theorem~\ref{relax-C5-thm} without changing the subgroups of the $\mathcal{G}(x^*)$ that arise as line stabilisers.

Theorem~\ref{relax-C5-thm} now furnishes an essentially free affine action of $\pi_1(\mathcal{G})$ on a $\Lambda$-tree. That the line stabilisers with respect to this action are soluble follows from Lemma~\ref{sol-line-stabs}. By \cite[Proposition 25(1)]{affine-paper}, non-trivial soluble subgroups of $\Gamma$ stabilise a unique line.\qed

\begin{corollary}\label{conj-pinched}(cf \cite[4.19]{Bass})\\
Let $\Gamma_0$ be a group equipped with an essentially free $\theta$-affine action on a $\Lambda_0$-tree
where $\Lambda_0$ is a regular \oag\ and $\im\theta$ is soluble. Let $s_i$ and $t_i$ ($i\in I$) be
elements of $\Gamma_0$ such that \be\item $s_i$ and $t_i$ generate maximal soluble subgroups of $\Gamma_0$ and $\theta_{s_i}=\theta_{t_i}=1$ for all $i$, \item $s_i$ is not conjugate to $t_i^{-1}$, \item for $i\neq j$ and $w=\pm 1$, we cannot have $s_i$ conjugate to $s_j^w$ (resp. $t_j^w$) and $t_i$ conjugate to $t_j^{-w}$ (resp. $s_j^{-w}$);
and \item no three groups of the form $\langle s_i\rangle$ or $\langle t_i\rangle$ are conjugate in $\Gamma_0$.
\ee
Then the iterated HNN extension $$\langle \Gamma_0,u_i (i\in I)\ |\ u_is_iu_i^{-1}=t_i (i\in I)\rangle$$
is $\ATF(\Z\times\Lambda_0)$.
\end{corollary}

\pf
Take a graph of groups $\mathcal{G}$ consisting of one vertex and with one edge $e_i$ for each element of $I$, with edge groups $\langle x_i\rangle$ and embedding maps given by $\alpha_{e_i}(x_i)=s_i$ and $\alpha_{\bar{e}_i}(x_i)=t_i$. In the notation of Theorem~\ref{Bass3.8'}, we take $g_{e_i}=u_i$. Then $\Gamma=\pi_1(\mathcal{G})$. By \cite[Corollary 26(1)]{affine-paper} and \cite[Proposition 27]{affine-paper}, the embedded images of edge groups coincide with end stabilisers with respect to the action of $\Gamma_0$.
Now all edge groups are cyclic, and condition (d) guarantees that (C2') is satisfied. Moreover, conditions (b), (c) and (d) ensure that no $\alpha_e(s)$ is conjugate to $\alpha_e(s^{-1})$; the details are similar to the argument given in \cite[4.19]{Bass}.
It follows from
Theorem~\ref{relax-C5-thm} that $\Gamma$ is $\ATF(\Z\times\Lambda_0)$. \qed

Theorem~\ref{intro-thm} of the introduction is a special case of Corollary~\ref{conj-pinched}.

Like Theorem~\ref{relax-C5-thm} the following result concerns combinations of ATF groups, but this time exploiting the degree of freedom afforded by the choice of $\mu$ to ensure that (C5) is satisfied.

\begin{theorem}\label{C5-mu}
Let $\Gamma$ be a group.
Consider the following conditions.
\be\item
$\Gamma$ admits a free affine action on a $\Z\times\Lambda_0$-tree.
\item There exists a graph of groups $(\mathcal{G},Y^*,T^*)$
such that $\Gamma\cong\pi_1(\mathcal{G})$ and the following are satisfied.
\be\item each vertex group $\mathcal{G}(x^*)$ is free on $\Rho_{x^*}=\{\rho_{{x^*},i}:i\in I_{x^*}\}$. Let $I=\coprod I_{x^*}$ and $\Rho=\coprod \Rho_{x^*}$, so that
$\Rho=\{\rho_i:i\in I\}$;
\item each vertex group $\mathcal{G}(x^*)$ has a free isometric action on a $\Lambda_0$-tree with hyperbolic length function $\ell=\ell^{x^*}$;
\item each $\alpha_e\mathcal{G}(e)$ is a maximal cyclic subgroup of $\mathcal{G}(\partial_0e)$, condition (C2') is satisfied, and no $\alpha_e(s)$ is conjugate in $\Gamma$ to $\alpha_e(s^{-1})$;
\item Let $E^+$ be an orientation of $E(Y^*)$, for $e\in E^+$ let $s=s_e=s_{\bar{e}}$ a generator of $\mathcal{G}(e)$ as in Lemma~\ref{edge-ends}, and put $u=u_e=\alpha_e(s_e)$ and $v=v_e=\alpha_{\bar{e}}(s_e)$.
For $x^*\in Y^*$, and $g\in F=\mathcal{G}(x^*)$, let $g_{i}$ ($i\in I$) be the integers such that $gF'=\prod_{i\in I}\rho_i^{g_i}F'$ where $F'=[F,F]$. (Of course $g_i=0$ for all but finitely many $i$).

There are elements $\bar{\rho}\in\Lambda_0$ for each ${\rho}\in \Rho$ satisfying the following.\\

\hspace{-1cm}\emph{(C5')}\h $\sum_{i\in I_{\partial_0e}}u_i\bar{\rho_i}=\sum_{i\in I_{\partial_0\bar{e}}}v_i\bar{\rho_i}=\ell^{\partial_0e}(v)-\ell^{\partial_0\bar{e}}(u)\h e\in E(Y^*)$

\ee \ee

Then (b)$\Rightarrow$(a). If $\Lambda_0=\Z$ then (a)$\Rightarrow$(b).

\end{theorem}

\pf
(a)$\Rightarrow$(b) assuming $\Lambda_0=\Z$: We apply Theorem~\ref{Bass3.5'}, noting that a $\Z$-free group is a free group, and taking $\theta_\gamma=1$ for all $\gamma\in\Gamma$. Thus (S1)-(S5) are given, satisfying (C1)-(C5). It suffices to show that (C5') is satisfied. Applying Lemma~\ref{edge-ends} we have an orientation $E^+\subseteq E(Y^*)$ and an assignment $e\mapsto\epsilon_e$ of ends such that for $e\in E^+$ with $\mathcal{G}(e)\neq 1$,
we have $\tau_{\epsilon_e}\alpha_e(s_e)>0$ and
$\tau_{\epsilon_{\bar{e}}}\alpha_{\bar{e}}(s_{e})<0$. Conditions (C4) and (C5) guaranteed by Theorem~\ref{Bass3.5'} give $-\left(\tau_e\alpha_e(s)+\tau_{\bar{e}}\alpha_{\bar{e}}(s)\right)=\mu_{\alpha_e(s)}=\mu_{\alpha_{\bar{e}(s)}}$. Moreover, $\mu_u=\mu_{\alpha_e(s)}$ and $\mu_v=\mu_{\alpha_{\bar{e}}(s)}$ and $\ell^{\partial_0e}(u)=\tau_e\alpha_e(s)$ and $\ell^{\partial_0\bar{e}}(v)=-\tau_{\bar{e}}\alpha_{\bar{e}}(s)$.
Next note that
$\mu:\Gamma\to\Z$ is a homomorphism. Put $\bar{\rho}=\mu_\rho$ for $\rho\in \Rho$. Then $\mu_\gamma=\sum_{i\in I}\gamma_i\bar{\rho}_i$ for $\gamma\in\Gamma$. Condition (C5') now follows.\\

(b)$\Rightarrow$(a):
We use Theorem~\ref{Bass3.8'}. Take all $\theta_{g_e}$ and $\theta^{x^*}$ to be trivial, and put $\mu_{g_e}=0$ for all $e$. (In fact $\mu_{g_e}$ can be chosen arbitrarily for $e\notin E(T^*)$.) The data (S1)-(S4) are given, and Lemma~\ref{edge-ends} gives a choice of ends $\epsilon_e$. Taking arbitrary end maps $\delta_e$ towards $\epsilon_e$ ($e\in E(Y^*)$) we have (S5) and (C2) satisfied.
Replacing the given $\Lambda_0$-trees by their respective fulfilments, we also have (C1).

Condition (C3) is trivial in the case at hand, and (C4) reduces to the requirement that $\mu_{\alpha_e(s)}=\mu_{\alpha_{\bar{e}}(s)}$ for $s\in \mathcal{G}(e)$ and that $\mu$ restricts to a homomorphism on each vertex group $\mathcal{G}(x^*)$. Such a homomorphism is determined by the images of $\rho_{x^*,i}$, so we set $\mu_{\rho}=\bar{\rho}$ for $\rho\in \Rho_{x^*}$ as in (iv), and extend $\mu$ to $\mathcal{G}(x^*)$ for each $x^*\in Y^*$. Condition (C5') ensures that $\mu_{\alpha_e(s)}=\mu_u=\mu_v=\mu_{\alpha_{\bar{e}}(s)}$. Thus (C4) is satisfied.

Condition (C5) reduces to $$-\left(\tau_e\alpha_e(s)+\tau_{\bar{e}}\alpha_{\bar{e}}(s)\right)=\mu_{\alpha_e(s)}.$$
Our choice of ends $\epsilon_e$ ensures that
for $e\in E^+$ the element $u=\alpha_e(s)$ translates towards $\epsilon_e$ and $v=\alpha_{\bar{e}}(s)$ translates away from $\epsilon_{\bar{e}}$, whence $\ell(u)=\tau_e(u)=\tau_e\alpha_e(s)$ and $\ell(v)=-\tau_{\bar{e}}(v)=-\tau_{\bar{e}}\alpha_{\bar{e}}(s)$. Condition (C5') now gives $\mu_{\alpha_e(s)}=\ell(v)-\ell(u)=-(\tau_e\alpha_e(s)+\tau_{\bar{e}}\alpha_{\bar{e}}(s))$, which implies (C5).
\qed

Note that in case $u$ is not conjugate to $v^{-1}$ but $u,v,\in F'$, Theorem~\ref{C5-mu} affords a free affine action of $\langle F,s\ |\ sus^{-1}=v\rangle$ on a
$(\Z\times\Z)$-tree if and only if $\ell(u)=\ell(v)$ in which case the action is isometric.

\section{Relative hyperbolicity for $\ATF^o(\Z^n)$ groups}\label{rel-hyper-section}

\begin{theorem}
Let $\Gamma$ be a finitely generated $\ATF^o(\Z^n)$ group. (See Notation~\ref{notation}.) \be\item $\Gamma$ is relatively hyperbolic with torsion-free nilpotent parabolic subgroups.
\item $\Gamma$ has solvable word problem.
\item $\Gamma$ has solvable conjugacy problem.
\item $\Gamma$ has solvable isomorphism problem.
\item $\Gamma$ is locally quasiconvex.\ee
\end{theorem}

\pf
(a) Induction on $n$. The assertion is trivial if $n=1$, so inductively suppose that finitely generated $\ATF(\Z^{n-1})$ groups are relatively hyperbolic with nilpotent parabolic subgroups.
Then $\Gamma$ admits a graph of groups decomposition as in Theorem~\ref{Bass3.5'}
where each vertex group is $\ATF(\Z^{n-1})$ and each edge group is finitely generated torsion free nilpotent. Since $\Gamma$ is finitely generated, the underlying graph is finite and the vertex groups are finitely generated and thus relatively hyperbolic. By \cite[Corollary 26(2)]{affine-paper}, all edge groups $\mathcal{G}(e)$ are either trivial, (proper) maximal nilpotent or bijectively mapped onto the end vertex group $\alpha_e\mathcal{G}(e)$. Inductively we will adjust the graph of groups decomposition given by the affine action by collapsing an edge at a time so that at each stage the vertex group is relatively hyperbolic with maximal nilpotent parabolics, and the fundamental group of the graph of groups obtained is equal to $\Gamma$.

For each $e$ for which $\mathcal{G}(e)=1$, we can replace the subgraph consisting of the edge $e$ and its endpoints by a single vertex $x^*$, with vertex group $$\mathcal{G}(x^*)=\left\{\begin{array}{rl}\mathcal{G}(\partial_0e)\ast\mathcal{G}(\partial_0\bar{e}) & \partial_0e\neq\partial_0\bar{e}\\ \mathcal{G}(\partial_0e)\ast\langle g_e\rangle & \partial_0e=\partial_0\bar{e}\end{array}\right.$$
Note that $\mathcal{G}(x^*)$, as a free product of relatively hyperbolic groups, is itself relatively hyperbolic with the same parabolic subgroups. Moreover for other edges $f$ with $\partial_0f=\partial_ie$ ($i=0,1$), the edge group embedding $\alpha_f$ can be adjusted to embed in $\mathcal{G}(x^*)$ in an obvious way. Since the non-cyclic nilpotent subgroups of a free product are contained in conjugates of a vertex group, if $\alpha_f\mathcal{G}(f)$ is maximal nilpotent in $\mathcal{G}(\partial_0f)$, the embedded image of $\mathcal{G}(f)$ in $\mathcal{G}(x^*)$ is also maximal nilpotent.

Now suppose that $\mathcal{G}(e)\neq 1$. If $\partial_0e\neq\partial_0\bar{e}$ and $\alpha_e\mathcal{G}(e)=\mathcal{G}(\partial_0e)$, then we can replace the subgraph consisting of the edge $e$ and its endpoints by a single vertex $x^*$, with vertex group $\mathcal{G}(x^*)=\mathcal{G}(\partial_0\bar{e})$, much as in the previous case considered. If $\partial_0e=\partial_0\bar{e}$ and $\alpha_e\mathcal{G}(e)=\mathcal{G}(\partial_0e)$, then
since $\alpha_{\bar{e}}\mathcal{G}({e})$ is maximal nilpotent in $\mathcal{G}(\partial_0e)$, we must have $\alpha_{\bar{e}}\mathcal{G}(e)=\mathcal{G}(\partial_0e)$. Now $\alpha_e\mathcal{G}(e)$ as a non-trivial nilpotent subgroup of $\Gamma$ stabilises a unique line of $X$. Since this subgroup is normalised by $g_e$, the line is also stabilised by the fundamental group $\Gamma_0=\langle\mathcal{G}(\partial_0e),g_e\ |\ g_e\alpha_e(s)g_e^{-1}=g_{\bar{e}}\alpha_{\bar{e}}(s)g_{\bar{e}}^{-1}\ s\in\mathcal{G}(e)\rangle$ of the subgraph spanned by the edge $e$. Thus $\Gamma_0$ is nilpotent by \cite[Proposition 25(3)]{affine-paper}. We adjust the peripheral structure of $\Gamma$ in this case by removing $\alpha_e\mathcal{G}(e)$ (if necessary) and adding the subgroup $\Gamma_0$.

Finally suppose that $\mathcal{G}(e)\neq 1$,
and $\alpha_e\mathcal{G}(e)$ is a proper maximal nilpotent subgroup of $\mathcal{G}(\partial_0e)$. Let $\Gamma_0$ be the fundamental group of the subgraph spanned by $e$. If $\partial_0e\neq\partial_0\bar{e}$ then by \cite[Corollary 1.7(1)]{Bigdely-Wise}, the fundamental group of the subgraph spanned by $e$ is relatively hyperbolic with maximal nilpotent parabolics. Otherwise $\partial_0e=\partial_0\bar{e}=x^*$, say. Let $\Gamma_0$ be the fundamental group of the subgraph of groups spanned by $e$. By \cite[Proposition 28(2)]{affine-paper} $\mathcal{G}(x^*)$ is CSN, and $\alpha_e\mathcal{G}(e)$
and $\alpha_{\bar{e}}\mathcal{G}(e)$ are proper maximal nilpotent in $\mathcal{G}(x^*)$. If these subgroups are not conjugate in $\mathcal{G}(x^*)$ it follows from \cite[Corollary 1.7(2)]{Bigdely-Wise} (or by \cite[Theorem 0.1(3')]{Dahmani-combination}) that $\Gamma_0$ is relatively hyperbolic, so suppose that $u\in \mathcal{G}(x^*)$ with $u\alpha_{\bar{e}}(s)u^{-1}=\alpha_e(s)$ for all $s\in\mathcal{G}(e)$. Then $\Gamma_0$ is presented by $$\langle \mathcal{G}(x^*), g_e\ |\ g_eu\alpha_{\bar{e}}(s)u^{-1}g_e^{-1}=\alpha_{\bar{e}}(s)\ (s\in\mathcal{G}(e)\rangle.$$
Replacing $g_e$ by $g_eu$ we lose no generality in assuming that $\alpha_e\mathcal{G}(e)=\alpha_{\bar{e}}\mathcal{G}(e)$. It follows from \cite[Corollary 1.7(3)]{Bigdely-Wise} that $\Gamma$ is relatively hyperbolic with nilpotent parabolic subgroups.
This proves part (a).

Parts (b), (c) and (d) follow from part (a) together with the fact that the respective problems are solvable in the class of relatively hyperbolic groups with nilpotent parabolics. The respective assertions required are \cite[Theorem 4.14]{Farb-rel-hyper}, \cite[Theorem 1.1]{Bumagin} and \cite[Theorem 1.5]{Dahmani-Touikan}.

Part (e) follows from part (a) and \cite[Theorem D]{Bigdely-Wise}, once one notes that nilpotent groups are Noetherian, and $\Gamma$ has a small hierarchy --- in other words, $\Gamma$ can be built up from groups containing no non-abelian free group by a finite sequence of amalgamated free products and HNN extensions along small subgroups.\qed\\

\section{A class of one-relator groups}\label{one-rel-section}

We consider here a class of groups that neatly illustrates a dividing line between $\ITF$ and $\ATF$ groups.

Denote the commutator $x^{-1}y^{-1}xy$ by $[x,y]$.

\begin{lemma}\label{lengths-comms}
Let $G$ be a group equipped with a free isometric action on a $\Lambda_0$-tree, let $x,y\in G$ and suppose that $xy\neq yx$. If $A_x\cap A_y$ contains at most one point, then $\ell(x^{-1}y^{-1}xy)=2\ell(x)+2\ell(y)+4d(A_x,A_y)$. (Here $d(A_x,A_y)=0$ if the axes are not disjoint.) If $A_x\cap A_y$ is a segment $[l,r]$ and $\Xi=d(l,r)>0$ then $\ell(x^{-1}y^{-1}xy)=2\ell(x)+2\ell(y)-2\Xi$.
\end{lemma}

\pf Suppose first that $A_x\cap A_y$ contains at most one point.
Then $y^{-1}\cdot (A_x\cap A_y)=A_{y^{-1}xy}\cap A_y$ contains at most
one point, and this point cannot also lie in $A_x\cap A_y$, since
$y$ is hyperbolic. Thus, by \cite[Lemma 2.1.11]{Chiswell-book} $A_x$ and $A_{y^{-1}xy}$ are disjoint, and
the closed bridge joining $A_x$ and $A_{y^{-1}xy}$ has the form
$[p,q,y^{-1}q,y^{-1}p]$ where $[p,q]$ is the bridge joining $A_x$ and $A_y$.
Therefore $d(A_x,A_{y^{-1}xy})=\ell(y)+2d(A_x,A_y)$.

Using \cite[Lemma 3.2.2]{Chiswell-book} we now obtain\\
\begin{eqnarray*}\ell(x^{-1}(y^{-1}xy))&=&\ell(x^{-1})+\ell(y^{-1}xy)+2d(A_x,A_{y^{-1}xy})\\
&=&2\ell(x)+2\ell(y)+4d(A_x,A_y).\end{eqnarray*}

Suppose now that $A_x\cap A_y$ is a closed segment $[l,r]$ with $\Xi=d(l,r)>0$. Since $\ell(x^{-1}y^{-1}xy)=\ell(y^{-1}x^{-1}yx)$, we lose no generality in assuming $\ell(y)\geq\ell(x)$. Further, since $xy\neq yx$ the intersection $[l,r]$ must have length strictly less than $\ell(x)+\ell(y)$, since otherwise $l$ or $r$ is fixed by some commutator of the form $[x^{\pm 1},y^{\pm 1}]$. Suppose next that $A_x\cap A_{y^{-1}xy}$ contains more than one point. Then the length of the intersection $A_x\cap A_{y^{-1}xy}$ is equal to $\Xi-\ell(y)$ and applying \cite[Lemma 3.3.3(2)]{Chiswell-book} to $g=x^{-1}$ and $h=(y^{-1}xy)^{-1}$, we obtain \begin{eqnarray*}\ell(x^{-1}y^{-1}xy)&=&\ell(x^{-1})+\ell(y^{-1}xy)-2(\Xi-\ell(y))\\ &=&2\ell(x)+2\ell(y)-2\Xi.\end{eqnarray*}

Finally if $A_x\cap A_{y^{-1}xy}$ contains at most one point then
$\Xi\leq\ell(y)$ and $d(A_x,A_{y^{-1}xy})=\ell(y)-\Xi$ so using \cite[Lemma 3.3.3(2)]{Chiswell-book} or \cite[Lemma 3.2.2]{Chiswell-book} we obtain
\begin{eqnarray*}\ell(x^{-1}y^{-1}xy)&=&\ell(x^{-1})+\ell(y^{-1}xy)+2(\ell(y)-\Xi)\\
&=&2\ell(x)+2\ell(y)-2\Xi.\end{eqnarray*}\qed

Note that Chiswell considers the hyperbolic length of a commutator
in a remark at the end of \cite[\S3.3]{Chiswell-book}; unfortunately
there is an error in the discussion there. (In the notation used
there, $gh$ and $g^{-1}h^{-1}$ meet coherently, not $gh$ and $hg$.)

\begin{corollary}\label{xmyn-length}
If $G$ has a free isometric action on a $\Lambda_0$-tree, and $x,y\in G$ with $xy\neq yx$ then $$\ell([x^m,y^n])=\left\{\begin{array}{ll}2|m|\ell(x)+2|n|\ell(y)+4d(A_x,A_y) & \mbox{if }A_x \mbox{ and }A_y\mbox{ have at most one point in common}\\
2|m|\ell(x)+2|n|\ell(y)-2\Xi & \mbox{otherwise.}\end{array}\right.$$\qed
\end{corollary}

Let $m,n,r$ and $s$ be non-zero integers and consider the group $$\Gamma(m,n;r,s)=\langle x,y,t\ |\ t[x^m,y^n]t^{-1}=[x^r,y^s]\rangle.$$

\begin{theorem}\label{xmyn-action}
Let $m,n,r$ and $s$ be non-zero integers.
\be
\item
The following assertions are equivalent.
\be
\item $\Gamma(m,n;r,s)$ is $\ATF^e$ (see Notation~\ref{notation}).
\item We do not have $m=-r$ and $n=s$, or $m=r$ and $n=-s$.
\item $\Gamma(m,n;r,s)$ is $\ATF^e(\Z\times\Q)$.
\item $\Gamma(m,n;r,s)$ is $\ATF^e(\Z^3)$.
\ee

\item
Assume that we do not have $m=-r$ and $n=s$, or $m=r$ and $n=-s$.
The following are equivalent.
\be\item $\Gamma(m,n;r,s)$ is $\ITF$.
\item $|m|-|r|$ and $|s|-|n|$ have the same sign.
\item $\Gamma(m,n;r,s)$ is $\ITF(\Z\times\Z)$.
\ee

\ee

\end{theorem}

\pf Write $\Gamma=\Gamma(m,n;r,s)$.

(a)
A violation of condition (ii) would imply that $[x^m,y^n]$ is conjugate in $\Gamma$ to its inverse, which is impossible in an $\ATF^e$ group, by Lemma~\ref{ess-hyper-props}(f).
Therefore (i)$\Rightarrow$(ii).

Now assume that condition (ii) holds.
Take any free isometric action of the free group $F$ on $\{x,y\}$ on a $\Q$-tree (for example the Base Change functor $X'=\Q\otimes_{\Z}X$ applied to the Cayley graph of $F$, viewed as a $\Z$-tree.) Condition (ii) is immediate from Theorem~\ref{relax-C5-thm} once one notes that $[x^m,y^n]$ is not conjugate in $F(x,y)$ to the inverse of $[x^r,y^s]$. Thus (ii)$\Rightarrow$(iii).

Again, assume (ii) and let $F_1$ be the free group on $\{x,z\}$. Starting from the Cayley graph with respect to this basis, assign a length of $(0,2)\in\Z^2$ to each edge labelled by $x$, and a length of $(1,0)$ to each edge labelled by $z$; this gives rise to a $\Z^2$-tree $Y$ on which $F_1$ has a natural free isometric action.
Put $y=z^{-1}xz$ and $F=\langle x,y\rangle$, and note that with
respect to the induced action of $F$ on $Y$, we have
$\ell(x)=(0,2)=\ell(y)$, the axes $A_x$ and $A_y$ are disjoint, and
$d(A_x,A_y)=(1,0)$. By Corollary~\ref{xmyn-length} we have
$\ell([x^m,y^n])=2|m|\ell(x)+2|n|\ell(y)+4d(A_x,A_y)=(4,4|m|+4|n|)$,
and $\ell([x^r,y^s])=(4,4|r|+4|s|)$. Put $a=
\left(|r|-|m|+|s|-|n|\right)$, and $\theta_{g_e}:(p,q)\mapsto(p,q+ap)$.
Now take the natural graph
of groups corresponding to the presentation of $\Gamma$: one vertex
$x^*$ with $\mathcal{G}(x^*)=F$, one edge $e$ with $\theta_{g_e}$ as just
described, $t=g_e$ as the element corresponding to the edge $e$, $w$
as a generator of the (infinite cyclic) edge group,
$\alpha_e(w)=[x^m,y^n]$ and $\alpha_{\bar{e}}(w)=[x^r,y^s]$. Let $\epsilon_e$ be the end of $Y$ towards which $[x^m,y^n]$ translates, and $\epsilon_{\bar{e}}$ the end of $Y$ towards which the inverse of $[x^r,y^s]$ translates.
Now
$\Gamma=\pi_1(\mathcal{G})$, and by Theorem~\ref{Bass3.8'}  $\Gamma$
has a free $\beta$-affine action on a $\Z^3$-tree with $\beta_g=1$
for $g$ in $F$ in $\Gamma$, and $\beta_t(1,\lambda_0)=(1,\theta_{g_e}(\lambda_0))$ ($\lambda_0\in\Z^2$).
This shows that (ii)$\Rightarrow$(iv).

It is clear that (iii) or (iv) imply (i).

(b)
Fix a free isometric action of $\Gamma=\Gamma(m,n;r,s)$ on a
$\Lambda_0$-tree. Since the hyperbolic lengths of conjugate
elements are equal, we must have $\ell([x^m,y^n])=\ell([x^r,y^s])$
giving $|m|\ell(x)+|n|\ell(y)=|r|\ell(x)+|s|\ell(y)$ by
Corollary~\ref{xmyn-length}. Thus
$(|m|-|r|)\ell(x)=(|s|-|n|)\ell(y)$. Since $\ell(x)$ and $\ell(y)$
are both positive, we must have either $|m|-|r|=|s|-|n|=0$, or
$|m|-|r|$ and $|s|-|n|$ both positive, or $|m|-|r|$ and $|s|-|n|$
both negative. This shows that (i)$\Rightarrow$(ii).

Assume (ii); suppose in particular that $|m|-|r|$ and $|s|-|n|$ are both positive (the
case where both are negative is easily reduced to this case). Note
that the commutators $[x^m,y^n]$ and $[x^r,y^s]$ can only be
conjugate in $F$ if $m=r$ and $n=s$, or $m=-r$ and $n=-s$, which is
impossible in the case at hand. Since these commutators are not
proper powers in $F(x,y)$, and not conjugate to each other, the
result \cite[4.19]{Bass} can be applied to obtain a free isometric
action of $\Gamma$ on a $\Z\times\Z$-tree: it suffices to show that
there is a free isometric action of the free group $F=F(x,y)$ on a
$\Z$-tree with $\ell([x^m,y^n])=\ell([x^r,y^s])$. So consider a free
action of $F$ on a $\Z$-tree with $\ell(x)=|s|-|n|$ and
$\ell(y)=|m|-|r|$; this can be done by taking the Cayley graph of
$F$ on $\{x,y\}$ and assigning the length $|s|-|n|$ to each edge
labelled $x$ and length $|m|-|r|$ to those labelled $y$. Corollary~\ref{xmyn-length} now
gives $\ell([x^m,y^n])=\ell([x^r,y^s])$.

If $m=r$ and $s=n$, then $\Gamma$ reduces to a benign HNN extension
of $F$ as in \cite[4.16]{Bass}. Thus $\Gamma$ has a free isometric
action as required. If $m=-r$ and $s=-n$ we have
$[x^m,y^n]=y^{-n}x^{-m}[x^r,y^s]x^{m}y^{n}$ so that
$[x^r,y^s]=t[x^m,y^n]t^{-1}$ is centralised by $x^{m}y^{n}t^{-1}$ in
$\Gamma$. Now replacing $t$ by $\bar{t}=x^{m}y^{n}t^{-1}$ in the
presentation of $\Gamma$, we obtain a benign HNN extension of $F$
which is isomorphic to $\Gamma$, giving, once again, a free
isometric action of $\Gamma$ on a $\Z\times\Z$-tree. Thus (ii)$\Rightarrow$(iii).

The equivalence of assertions in part (b) is now clear.
\qed

It is shown in \cite[Theorem 1(3)]{affine-paper} that the Heisenberg group
$\UT(3,\Z)$ is an example of an $\ATF(\Z^3)$ group that is not
$\ITF$ and that for $m>1$ the soluble
Baumslag-Solitar groups $\BS(1,m)$ are $\ATF(\Z\times\R)$
but not $\ATF(\Z^n)$ for any $n$, nor $\ITF$.

We conclude this section with another interesting example of an $\ITF$
group. Although not directly relevant to affine actions,
it uses ideas discussed in this section.

\begin{theorem}
The group $$\Gamma_1=\langle x,y,t\ |\ txt^{-1}=[x,y]\rangle.$$
is $\ITF(\Z^2)$, virtually special (and hence
virtually residually torsion-free nilpotent), but not residually
nilpotent.
\end{theorem}

\pf
Fix the standard free isometric action of the free group
$F_0=F(u,v)$ on its Cayley graph (with respect to the basis
$\{u,v\}$), and put $x_1=u^3v$ and $y_1=u$. Then
$\ell(x_1)=3\ell(u)+\ell(v)=4$ and $\ell(y_1)=1$, while
$\Xi(x_1,y_1)=\ell(u^3)=3$. Thus by Lemma~\ref{lengths-comms}, we
have
\begin{eqnarray*}\ell[x_1,y_1]&=&2\ell(x_1)+2\ell(y_1)-2\Xi(x_1,y_1)\\
&=&4\\ &=&\ell(x_1).\end{eqnarray*} Since $\{x_1,y_1\}$ is a
generating set, and thus a basis of $F_0$, this shows that the free
group $F$ on $\{x,y\}$ admits a free isometric action on a $\Z$-tree
with $\ell(x)=\ell[x,y]$. Moreover, since $[x,y]$ and $x$ are not
proper powers in $F$, the HNN extension $\Gamma_1=\langle F,t\ |\
txt^{-1}=[x,y]\rangle$ admits a free isometric action on a
$\Z\times\Z$-tree, by~\cite[4.19]{Bass}. That finitely generated
$\ITF(\Z^n)$ groups are
virtually special follows from \cite[Corollary 19]{KMS-long-survey}.

We now show that $\Gamma_1$ is not residually nilpotent. Let $N$ be
a nilpotent group of class $\leq c$, and $\phi:\Gamma_1\to N$ a
homomorphism. Let $u_k=t^kyt^{-k}$, let $v_0=x$, and inductively put
$v_k=[v_{k-1},u_{k-1}]$. Then $v_k$ is a simple commutator of weight $k$
for all $k\geq 1$ and so $\phi(v_{c+1})=1$. We claim that
$v_k=t^kxt^{-k}$ for all $k\geq 1$.

For inductively, if $k\geq 2$ then
\begin{eqnarray*}v_k&=&[v_{k-1},u_{k-1}]\\
&=&[t^{k-1}xt^{-k+1},t^{k-1}yt^{-k+1}]\\ &=&t^{k-1}[x,y]t^{-k+1}\\
&=& t^kxt^{-k}.
\end{eqnarray*}
Thus $1=\phi(v_{c+1})=\phi(t^{c+1}xt^{-c-1})$, forcing $\phi(x)=1$. Thus the
image of $x$ under every homomorphism $\Gamma_1\to N$ is trivial.\qed

\bibliographystyle{plain}
\def\cprime{$'$}

\begin{minipage}[t]{3 in}
\noindent Shane O Rourke\\ Department of Mathematics\\
Cork Institute of Technology\\
Rossa Avenue\\ Cork\\ IRELAND
\\ \verb"shane.orourke@cit.ie"
\end{minipage}

\end{document}